\documentclass[11pt,reqno]{amsart}
\usepackage{amsfonts}
\usepackage{amsaddr}
\usepackage{graphicx}
\usepackage{amscd}
\usepackage{amsthm}
\usepackage{amstext}
\usepackage{amsmath}
\usepackage{float}
\usepackage{amssymb,latexsym}
\usepackage{epsfig}
\usepackage{mathrsfs}
\usepackage{subfigure}
\usepackage{subfig}
\usepackage{cite}
\usepackage{multicol}
\usepackage[margin=3cm]{geometry}
\DeclareGraphicsExtensions{.png,.pdf,.jpg,.gif}

\newtheorem{teo}{Theorem}[section]

\newtheorem{lema}{Lemma}[section]

\newtheorem{prop}{Proposition}[section]
\theoremstyle{definition}

\newtheorem{remark}{Remark}[section]

\theoremstyle{plain}

\numberwithin{equation}{section}

\begin{document}

\title[Bifurcation and Periodic solutions] {Bifurcation and periodic solutions to population models with two dependent delays}

\author[A. G\'omez]{Adri\'an G\'omez \textsuperscript{$\MakeLowercase{a}$}, Jos\'e Oyarce \textsuperscript{$\MakeLowercase{a,*}$}}
\email{agomez@ubiobio.cl}

\author[J. Oyarce]{}
\email{jooyarce@egresados.ubiobio.cl}

\address{$^a$Departamento de Matem\'atica, Facultad de Ciencias, Universidad del B\'io-B\'io, Casilla 5-C, Concepci\'on, VIII-Regi\'on, Chile}

\thanks{*Corresponding author}



\begin{abstract}
We investigate the scalar autonomous equation with two discrete delays 
$$
\dot{x}(t)=f(x(t),x(t-r),x(t-\sigma)),
$$
where $f:\mathbb{R}^3\rightarrow \mathbb{R}$ is a continuously differentiable non-linear function such that $f(0,0,0)=0$. It is shown that if the difference between the delays is constant, then one of the delays becomes a Hopf-bifurcation parameter and, in addition, the absolute stability of the trivial solution can be established. Moreover, the direction of the Hopf bifurcation and the stability of the bifurcating periodic solutions are determined by using normal form theory. The main results are applied to guarantee the existence of positive periodic solutions to Nicholson's blowflies and Mackey-Glass models, both with a delayed harvesting term. The conclusions are illustrated by numerical simulations. 
\end{abstract}

\maketitle 

\section{Introduction}\label{S1}

To describe periodic oscillations in the experiments  of Nicholson \cite{Nicholson} with the Australian sheep blowfly \emph{Lucilia cuprina}, Gurney \emph{et al.} \cite{Gurney} proposed the following model
\begin{equation}\label{Intro1}
\dot{x}(t)=-\delta x(t)+Px(t-r)e^{-x(t-r)}.
\end{equation}
Here $x(t)$ represents the population, $\delta$ is the mortality rate, $P$ is the maximum per capita daily egg production, and $r$ is the time taken from birth to maturity. In \cite{Wei2005}, the authors studied the model \eqref{Intro1} by taking the delay as a parameter and, consequently, a result dealing with the existence of periodic solutions to \eqref{Intro1} is presented in the paper (see \cite[Thm. 2.3]{Wei2005}). 

One of the open problems formulated by Berezansky \emph{et al.} \cite{braverman} is to investigate the Nicholson's blowflies model with a delayed linear harvesting term
\begin{equation}\label{Nicholson}
\dot{x}(t) = -\delta x(t)+Px(t-r)e^{-x(t-r)}-H x(t-\sigma),
\end{equation}
where the harvesting $Hx(t-\sigma)$ is a function of the delayed estimate of the true population. The results obtained in this paper are applied to find the sufficient  conditions for the existence of bifurcating periodic solutions to \eqref{Nicholson} (see Section \ref{S4}). Some general cases of \eqref{Intro1} can be found, e.g., in \cite{Amster,bravd,braverman} and the references therein. 

To describe the dynamical behaviour of red blood cells production, Mackey and Glass \cite{Mackey1} formulated the following single-delayed model
\begin{equation}\label{Intro2}
\dot{u}(t)=-\delta u(t)+ \frac{\beta\theta^n}{\theta^n+u^n(t-r)},
\end{equation}
where $u(t)$ represents the circulating density of red blood cells, $\delta$ is the loss rate of red cells, $\beta$ is the maximal red blood cell production rate that the body can approach at low circulating red blood cell numbers, $r$ is a maturation delay, $n$ is a positive exponent, and $\theta$ is a shape parameter. Let $u(t)=\theta x(t)$, then \eqref{Intro2} becomes
\begin{equation}\label{Intro3}
\dot{x}(t)=-\delta x(t) + \frac{P}{1+x^n(t-r)},
\end{equation}
where $P=\beta/\theta$. Similar to \cite{Wei2005}, in \cite{Wei2007} the authors investigated the existence of Hopf bifurcations at a positive equilibrium $x=x_*$ of \eqref{Intro3} by using the delay as a parameter (see \cite[Thm. 2.3]{Wei2007}). Some results about oscillation and global attractivity of solutions to \eqref{Intro2} can be found, e.g., in \cite{Gopalsamy2,Gyori}. The results obtained in this paper are also applied to investigate \eqref{Intro3} with a delayed harvesting strategy (see Section \ref{S4}). 

It is well-known that scalar differential equations with two delays are more realistic than single-delayed models and, in particular, they have significant importance in biological applications. We refer to \cite{Wei2021,Huang2019,berbrav2,Braddock,Gopalsamy,RuanWeiLi} for population models with two delays; to \cite{Belair,Lainscsek1,Lainscsek2} for neurological models with two delays; and to \cite{ruan} for a compound optical resonator with two delays. Therefore, the purpose of this paper is to investigate the following scalar autonomous equation 
\begin{equation}\label{ecfinal}
\dot{x}(t)=f(x(t),x(t-r),x(t-\sigma)), \ \ r\geq 0, \ \ \sigma \geq 0.
\end{equation}
Here, for $f(u_1, u_2, u_3)$, we assume that $f: \mathbb{R}^3\rightarrow \mathbb{R}$ is a continuously differentiable non-linear function such that  $f(0,0,0)=0$ and  
\begin{equation*}
-a= \frac{\partial f}{\partial u_1}(0,0,0), \ \ -b=\frac{\partial f}{\partial u_2}(0,0,0), \ \ -c= \frac{\partial f}{\partial u_3}(0,0,0),
\end{equation*}
where $a, b$ and $c$ are real numbers. We apply linear methods to study the local behaviour of the trivial solution to \eqref{ecfinal}, which leads us to investigate the linear equation 
\begin{equation}\label{Lineal}
\dot{y}(t)=- ay(t)-by(t-r)-cy(t-\sigma ).
\end{equation}
From a biological point of view, it is interesting to study positive periodic orbits that oscillate about positive equilibria of the delayed model. Obviously, this analysis can be done by shifting the positive equilibrium to zero and investigating a differential equation of the form \eqref{ecfinal} and its linearization \eqref{Lineal}.

In the last decades the delayed equation  \eqref{Lineal} has become of interest to many authors, e.g., Hale and Huang \cite{halegeometric} investigated the stability of the zero equilibrium of \eqref{ecfinal} and determined a geometrical description of the stable regions of \eqref{Lineal} in the $(r, \sigma)$--plane for some sets of $a, b$ and $c$. In the case $a=0, b>0,$ and $c>0$, Li  \emph{et al.}  \cite{RuanWeiLi} studied the local stability of the zero equilibrium of \eqref{ecfinal} and the existence of Hopf bifurcations considering one of the delays as a parameter, moreover the authors studied the direction of the bifurcation and the stability of the bifurcating periodic solutions by using the method of normal forms developed by Hassard  \emph{et al.}  \cite{Hassard}. Piotrowska \cite{Piotrowska} presented some remarks and corrections about the results in \cite{RuanWeiLi} and, in addition, the cases of \eqref{Lineal} with $a=0, b<0,$ and $ c<0$ or $a=0$ and $b c<0$ were  investigated. Note that in \cite{RuanWeiLi, Piotrowska} the bifurcation analysis is treated for two independent delays. Nevertheless, as we will see throughout this paper, this is not our case since we will assume that the difference between the delays is constant. If $bc=0, \ r \sigma=0$ or $r=\sigma$, then \eqref{Lineal} is equal to an equation with a single delay, which corresponds to already known results (see, e.g., \cite{Hale1, Hsmith}) and, therefore, in this paper we do not consider this latter cases. 

The main purpose of this paper is to study the existence and stability of periodic solutions to \eqref{ecfinal} when the parameters $a, b$ and $c$ belong to  suitable sets and the difference between the delays is constant, namely $\tau=r-\sigma$ with $\tau \in \mathbb{R}$ fixed. Based on some of the ideas and results in \cite{Hsmith,halegeometric, RuanWeiLi, ruan, Wei2005, Gu, Qi}, we introduce a new method to analyse the distribution of the roots of the characteristic equation corresponding to \eqref{Lineal} and, consequently, we prove the existence of local Hopf bifurcations about the zero equilibrium of \eqref{ecfinal}. Furthermore, by using the method of normal forms developed by Faria and Magalh\~aes \cite{Faria}, we study the direction of the Hopf bifurcation and the stability of the bifurcating periodic solutions. The results obtained in this paper are applied to Nicholson's blowflies and Mackey-Glass equations. 

The paper is structured as follows. In Section \ref{S2}, by choosing a delay as a parameter, we obtain the sufficient conditions to prove the existence of local Hopf bifurcations to \eqref{ecfinal} when the difference between the delays is constant and, in addition, an absolute stability criterion of the zero equilibrium is stated. In Section \ref{S3}, we study the direction of the bifurcation and the stability of the bifurcating periodic solutions by using normal forms. The results of Sections \ref{S2} and \ref{S3} are applied in Section \ref{S4} to prove the existence of positive periodic solutions to the Nicholson's blowflies population model \eqref{Nicholson}, and the Mackey-Glass model \eqref{Intro3} with a delayed harvesting term. A numerical analysis of the equation \eqref{Nicholson} is also presented in Section \ref{S4} using the Matlab dde23 Package. 

\section{Stability of the zero equilibrium and local Hopf bifurcations}\label{S2}

In this section, we prove main results about the existence of absolute stability and local Hopf bifurcations for the trivial solution of \eqref{ecfinal}.

We start this section by analysing the characteristic equation associated with \eqref{Lineal}, it is 
\begin{equation}\label{caracteristica}
h(\lambda) \stackrel{def}{=} \lambda + a + b e^{-\lambda r} + c e^{-\lambda \sigma}=0.
\end{equation}
Letting $\lambda = \mu + i\omega$ with $\omega\neq 0$, and  separating real and imaginary parts we have
\begin{gather}
\label{sistemacaracteristica}
\begin{aligned}
\mu +a &= - b e^{-\mu r}\cos (\omega r ) -c e^{-\mu \sigma }\cos ( \omega\sigma ), \\
\omega &=b e^{-\mu r} \sin (\omega r )+c e^{- \mu\sigma } \sin ( \omega \sigma) .
\end{aligned}
\end{gather}
We are interested in the existence of purely imaginary roots of \eqref{caracteristica}, hence we do  $\lambda = i\omega$ in \eqref{sistemacaracteristica} obtaining 
\begin{gather}\label{sistemaimaginariapura}
\begin{aligned}
a &=-b\cos (\omega r)  - c\cos (\omega \sigma ),  \\
\omega &= \ \  b \sin (\omega r) + c \sin (\omega \sigma ) . 
\end{aligned}
\end{gather}
If $a+b+c=0$, then $\omega=0$ is a solution of \eqref{sistemaimaginariapura} for all delays $(r,\sigma)\in \mathbb{R}_+^2$ and, therefore, we do not consider that case. Let
\begin{equation*}
X_+ \stackrel{def}{=} \lbrace (a,b,c)\in \mathbb{R}^3: a+b+c>0 \rbrace.
\end{equation*}
If $r=\sigma=0$ and $(a,b,c)\in X_+$, then the zero equilibrium of \eqref{ecfinal} is asymptotically stable. Also, from \cite[Prop. 2.1]{halegeometric}, if $(a,b,c)\in X_-=\lbrace (a,b,c)\in \mathbb{R}^3: a+b+c<0 \rbrace$, then the zero equilibrium of \eqref{ecfinal} is unstable for all delays $r\geq 0, \sigma \geq 0$. Thus, throughout this paper we assume 
\begin{description}
\item[(A)] $(a,b,c)\in X_+$,
\item[(B)] $bc\neq 0$. 
\end{description}

In order to find and to discard imaginary roots of (\ref{caracteristica}), we have the following lemma.
\begin{lema}\label{equivalentes} Let $\tau\stackrel{def}{=} r-\sigma \in \mathbb{R}$. 
\begin{enumerate}
\item[(i)] The system \eqref{sistemaimaginariapura} has a  root $\omega^*>0$ if and only if the equation 
\begin{equation}\label{equiv}
\cos(\omega \tau)=\dfrac{\omega^2}{2bc}+\dfrac{a^2-b^2-c^2}{2bc}
\end{equation}
has the same  root $\omega^*$.
\item[(ii)] A simple root $\omega^*>0$ of \eqref{equiv} is a simple root of \eqref{sistemaimaginariapura}.
\end{enumerate} 
\end{lema}

\begin{proof}(i) 
Adding up to square both sides of \eqref{sistemaimaginariapura} we arrive to \eqref{equiv}, hence it is clear that any root of the system \eqref{sistemaimaginariapura} is a root of equation \eqref{equiv}. Reciprocally, let $\omega>0$ be a root of  \eqref{equiv}. Rewriting \eqref{equiv} as

\begin{equation}\label{determinante}
2bc\cos(\omega \tau)+c^2+b^2=\omega^2+a^2,
\end{equation}
 we observe that $\Delta\stackrel{def}{=}-(2bc\cos(\omega \tau)+c^2+b^2)<0$. In consequence,  the system 

\begin{equation}\label{msystem}
\left(\begin{array}{ll}
-(c+b\cos(\omega\tau))&  b\sin(\omega\tau)\\  b\sin(\omega\tau)&(c+b\cos(\omega\tau))
\end{array}\right)\left(\begin{array}{l}
u_1\\ u_2
\end{array}\right)=\left(\begin{array}{l}
a\\ \omega
\end{array}\right)
\end{equation}
has a unique solution such that 
\begin{eqnarray*}
\left(\begin{array}{l}
u_1\\ u_2
\end{array}\right)&=&\dfrac{1}{|\Delta |}\left(\begin{array}{ll}
-(c+b\cos(\omega\tau))&  b\sin(\omega\tau)\\  b\sin(\omega\tau)&(c+b\cos(\omega\tau))
\end{array}\right)\left(\begin{array}{l}
a\\ \omega
\end{array}\right),\\
&=&\dfrac{1}{\sqrt{|\Delta|}}\left(\begin{array}{ll}
-(c+b\cos(\omega\tau))&  b\sin(\omega\tau)\\  b\sin(\omega\tau)&(c+b\cos(\omega\tau))
\end{array}\right)\left(\begin{array}{l}
\frac{a}{\sqrt{|\Delta|}}\\ \frac{\omega}{\sqrt{|\Delta|}}
\end{array}\right).
\end{eqnarray*}
Note that 
\begin{equation*}
det\left(\dfrac{1}{\sqrt{|\Delta|}}\left(\begin{array}{ll}
-(c+b\cos(\omega\tau))&  b\sin(\omega\tau)\\  b\sin(\omega\tau)&(c+b\cos(\omega\tau))
\end{array}\right)\right)=-1,
\end{equation*}
and from \eqref{determinante} it follows that
\begin{equation*}
\left\|\left(\begin{array}{l}\frac{a}{\sqrt{|\Delta|}}\\ \frac{\omega}{\sqrt{|\Delta|}}\end{array}\right)\right\|= \dfrac{1}{\sqrt{|\Delta |}}\sqrt{a^2+\omega^2}=1,
\end{equation*}
whence we obtain that $u_1^2+u_2^2=1$. Therefore, there exists $\sigma\in [0,2\pi/\omega)$ such that $u_1=\cos(\omega \sigma)$, $u_2=\sin(\omega \sigma)$ and,  consequently  

\begin{equation*}
\left(\begin{array}{ll}
-(c+b\cos(\omega\tau))&  b\sin(\omega\tau)\\  b\sin(\omega\tau)&(c+b\cos(\omega\tau))
\end{array}\right)\left(\begin{array}{l}
\cos(\omega\sigma)\\ \sin (\omega\sigma)
\end{array}\right)=\left(\begin{array}{l}
a\\ \omega
\end{array}\right).
\end{equation*}
The last system is equivalent to \eqref{sistemaimaginariapura}, thus we have shown  that $\omega$ solves \eqref{sistemaimaginariapura}.

In order to prove (ii), we assume now that $\omega $ is a multiple root of \eqref{sistemaimaginariapura}. Let $f_1(\omega) \stackrel{def}{=} -b\cos (\omega r)  - c\cos (\omega \sigma )$ and $ f_2(\omega)  \stackrel{def}{=} \ \  b \sin (\omega r) + c \sin (\omega \sigma ) $, then the system \eqref{sistemaimaginariapura} takes the form
\begin{equation}\label{sis1}
\begin{array}{l}a=f_1(\omega),\\
\omega=f_2(\omega),
\end{array}
\end{equation}
and \eqref{equiv} is
\begin{equation}\label{equ1}
a^2+\omega^2=f_1^2(\omega)+f_2^2(\omega).
\end{equation}
If $\omega$ is a multiple root of \eqref{sis1} it also satisfies $0=f'_1(\omega)$ and  $1=f'_2(\omega)$, hence
\begin{equation*}
\omega=f_1(\omega)f'_1(\omega)+f_2(\omega)f'_2(\omega),
\end{equation*}
which implies that $\omega$ is a multiple root of \eqref{equ1} and, therefore, a multiple root of \eqref{equiv}, showing (ii).
\end{proof}
The last lemma allows to find simple roots of \eqref{sistemaimaginariapura} by studying \eqref{equiv}. In order to find these roots, we enunciate the following.
 
\begin{lema}\label{lematau}
For  \eqref{sistemaimaginariapura} we have 
\begin{enumerate}
\item[(i)] If $bc>0$ and $|b-c|\leq |a|<|b+c|$, then for any $\tau\in\mathbb{R}$ there are smaller values of $ \omega^*>0$, $k_\tau\in \mathbb{N} _0$, and  $\bar{\sigma}\in [0,2\pi/\omega^*)$ such that for any $k\geq k_\tau$, we have delays 

 \begin{equation*}
 \begin{array}{l}
\sigma_k \stackrel{def}{=}\bar{\sigma}+\dfrac{2k\pi}{\omega^*}> 0,\\
\\
r_k \stackrel{def}{=}\tau+\sigma_k> 0,
 \end{array}
\end{equation*}  
such that $(r_k,  \sigma_k,\omega^*)$ solves  \eqref{sistemaimaginariapura}. Also, for each $(r_k,\sigma_k)$, $\omega^*$ is a simple root.
\item[(ii)] If $bc>0$ and $ |a|<|b-c|$, then for all $$\displaystyle|\tau|\in \left[0,\dfrac{\sqrt{(b-c)^2-a^2}}{bc}\right)\cup\bigcup_{n=0}^{\infty}\left[\dfrac{2n\pi}{\sqrt{(b-c)^2-a^2}},\dfrac{(2n+1)\pi}{\sqrt{(b-c)^2-a^2}}\right],$$ the conclusion in (i) is valid.
\item[(iii)] If $bc<0$, $|b+c|<|a|<|b-c|$,  there is a value $\tau^*> 0$ such that for $|\tau|< \tau^*$ the system \eqref{sistemaimaginariapura} has not any root $\omega>0$. And if  $|\tau|>\tau^*$, the conclusion in (i) is valid.
\item[(iv)] If $bc<0$ and $|a|<|b+c|$, then for all 

$$|\tau|\in  \left[0,\dfrac{\sqrt{(b+c)^2-a^2}}{|bc|}\right)\cup\bigcup_{k=0}^\infty \left[\dfrac{(2k-1)\pi}{\sqrt{(b+c)^2-a^2}},\dfrac{2k\pi}{\sqrt{(b+c)^2-a^2}}\right],$$
the conclusion in (i) is valid.
\end{enumerate}
\end{lema}  
\begin{proof}
Let $g(\omega) \stackrel{def}{=} \dfrac{\omega^2}{2bc}+\dfrac{a^2-b^2-c^2}{2bc}$, assuming $bc>0$ and $|b-c|\leq |a|<|b+c|$ we have that $g$ is an strictly increasing function for $\omega\in[0,\infty)$ such that $-1\leq g(0)<1$ and thus, $g(\omega^*)=\cos(\omega^* \tau)$ for a first positive $\omega^*$ and in this intersection the slopes have opposite signs, showing the simplicity of $\omega^*$. Now, from Lemma \ref{equivalentes}, if $\omega^*>0$ solves \eqref{equiv}, then $\Delta< 0$ and we can rewrite \eqref{msystem} as
 \begin{equation}\label{uvectores}
 \left(\begin{array}{l}
\cos(\omega^*\sigma)\\ \sin(\omega^*\sigma)
\end{array}\right)=\left(\begin{array}{l}\dfrac{b\omega^*\sin(\omega^*\tau)-a(c+b\cos(\omega^*\tau))}{{\omega^*} ^2+a^2}\vspace{0.3cm}\\

\dfrac{ab\sin(\omega^*\tau)+\omega^*(c+b\cos(\omega^*\tau))}{{\omega^*} ^2+a^2}
\end{array}\right).
 \end{equation}
As both sides of the last equation are unitary vectors, there is a first $\bar{\sigma}\in [0,2\pi/\omega^*)$ such that \eqref{uvectores} holds. Hence, by the periodicity of the left side of \eqref{uvectores}, any $\sigma_k=\bar{\sigma}+2k\pi/\omega^*$ with $k\in\mathbb{N}_0$, satisfies it. Therefore we can choose $k_\tau$ large enough such that $r_k\stackrel{def}{=}\tau+\sigma_k\geq 0$ for all $k\geq k_\tau$, and (i) holds.

Assuming $bc>0$ and $|a|<|b-c|$, then $g(0)<-1$ and $g(\omega)$ takes all the values 
in $[-1,1]$ when $\omega\in [\omega_1,\omega_2]$  with $g(\omega_{1,2})=\mp 1$. Hence, by continuity, the equation \eqref{equiv} has for each $\tau\in \mathbb{R}$ at least one root $\omega\in [\omega_1,\omega_2]$. Fixing $\tau$ and denoting the first root by $\omega^*$, in order to have a simple first root $\omega^*$ for \eqref{equiv}, we need to avoid 
\begin{equation}\label{tangente}
\dfrac{\omega^*}{bc}=-\tau\sin(\omega^* \tau).
\end{equation} 
Since $\omega^*\geq \omega_1$, if \eqref{tangente} is true, then 
 \begin{eqnarray*}
\dfrac{\omega_1}{bc}&= & \dfrac{\omega_1}{|bc|}
                  \leq \dfrac{\omega^*}{|bc|}
                  =|\tau| |\sin( \omega^*\tau)|
                  \leq |\tau|,
\end{eqnarray*}
which is imposible if $|\tau|<  \dfrac{\omega_1}{bc}=\dfrac{\sqrt{(b-c)^2-a^2}}{bc}$.
Also, if $|\tau|\neq 0$ satisfies 
\begin{equation}\label{posiciontau}
\dfrac{2k\pi}{|\tau|}\leq \omega_1<\dfrac{(2k+1)\pi}{|\tau|},
\end{equation}
then  $\omega^*\in\Big[\omega_1, \dfrac{(2k+1)\pi}{|\tau|}\Big)$ and $-\tau\sin(\omega^* \tau)<0<g(\omega^*)$. Consequently, from \eqref{posiciontau} if 
$
|\tau| \in \bigcup_{k=0}^\infty\Big[\dfrac{2k\pi}{\omega_1},\dfrac{(2k+1)\pi}{\omega_1}\Big)
$, 
then $\omega^*$ is a simple root and the conclusion in (i) is valid following the same proof.

Assuming $bc<0$ and $|b+c|<|a|<|b-c|$, then the function $g$ is strictly decreasing on $[0,\infty)$ and $g(0)\in (-1,1)$, hence the conclusion follows easily. 

Assuming $bc<0$ and $|a|<|b+c|$, we have that $g$ is strictly decreasing on $[0,\infty)$, and $g(0)>1$. By continuity, clearly \eqref{equiv} has at least one solution in $[\omega_1,\omega_2]$ with $g(\omega_{1,2})=\pm 1$. Then, we consider $\omega^*$ to be the first root of \eqref{equiv}.  In order to guarantee the simplicity of this first root,  we need to avoid \eqref{tangente} as above. In this case we take $$|\tau|<\dfrac{\omega_1}{|bc|}=\dfrac{\sqrt{(b+c)^2-a^2}}{|bc|}.$$
Then, $|\tau\sin(\tau\omega)|<\dfrac{\omega_1}{|bc|}\leq \dfrac{\omega}{|bc|}$ for all $\omega\in [\omega_1,\omega_2]$, and \eqref{tangente} is impossible. Also, if $\omega_1\in\left[\dfrac{(2n-1)\pi}{|\tau|},\dfrac{2n\pi}{|\tau|}\right]$ for some $n\in\mathbb{N}$, then as in the first case, the slopes at the point $\omega^*$ for the left and the right sides of \eqref{equiv} have opposite signs, hence $\omega^*$ is a simple root if 
\begin{equation*}
|\tau|\in \left[0,\dfrac{\sqrt{(b+c)^2-a^2}}{|bc|}\right)\cup\bigcup_{k=0}^\infty \left[\dfrac{(2k-1)\pi}{\sqrt{(b+c)^2-a^2}},\dfrac{2k\pi}{\sqrt{(b+c)^2-a^2}}\right],
\end{equation*} 
and the proof follows as in the case (i).
\end{proof}

If the parameters $a,b,c$ satisfy either \textbf{(A)}, $bc>0$ and $|a|> |b+c|$ or  \textbf{(A)},
 $bc<0$ and $|a|\geq |b-c|$, due to the non-existence of non-zero imaginary roots in \eqref{caracteristica}, it can be easily verified the absolute local stability of the trivial solution to \eqref{ecfinal}. 
 
\begin{figure}
{\centering
 \subfigure[$bc>0$ and $|b-c|\leq |a|<|b+c|$]{\scalebox{0.43}{\includegraphics{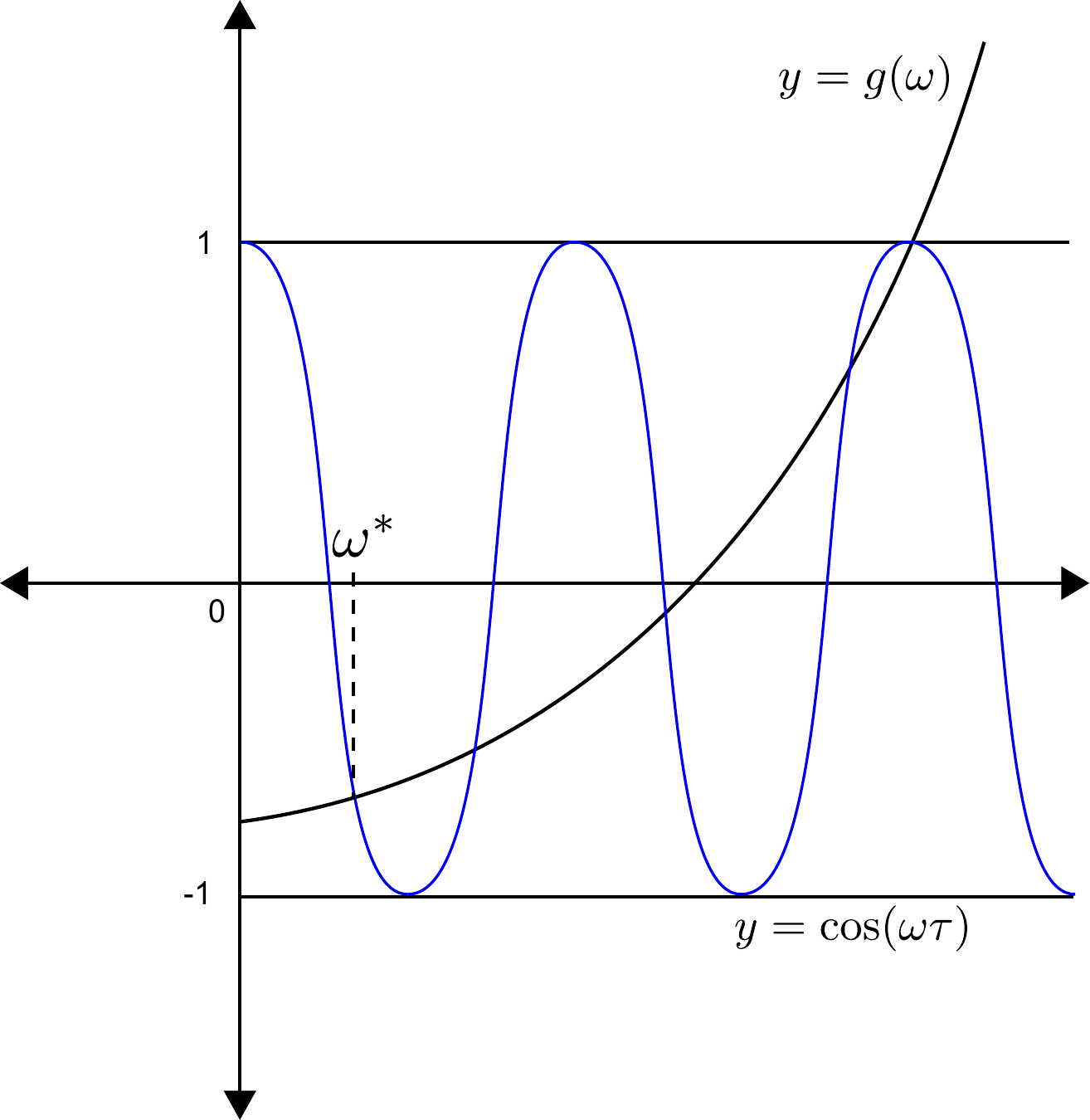}}}
       \subfigure[$bc>0$ and $|a|<|b-c|$]{\scalebox{0.43}{\includegraphics{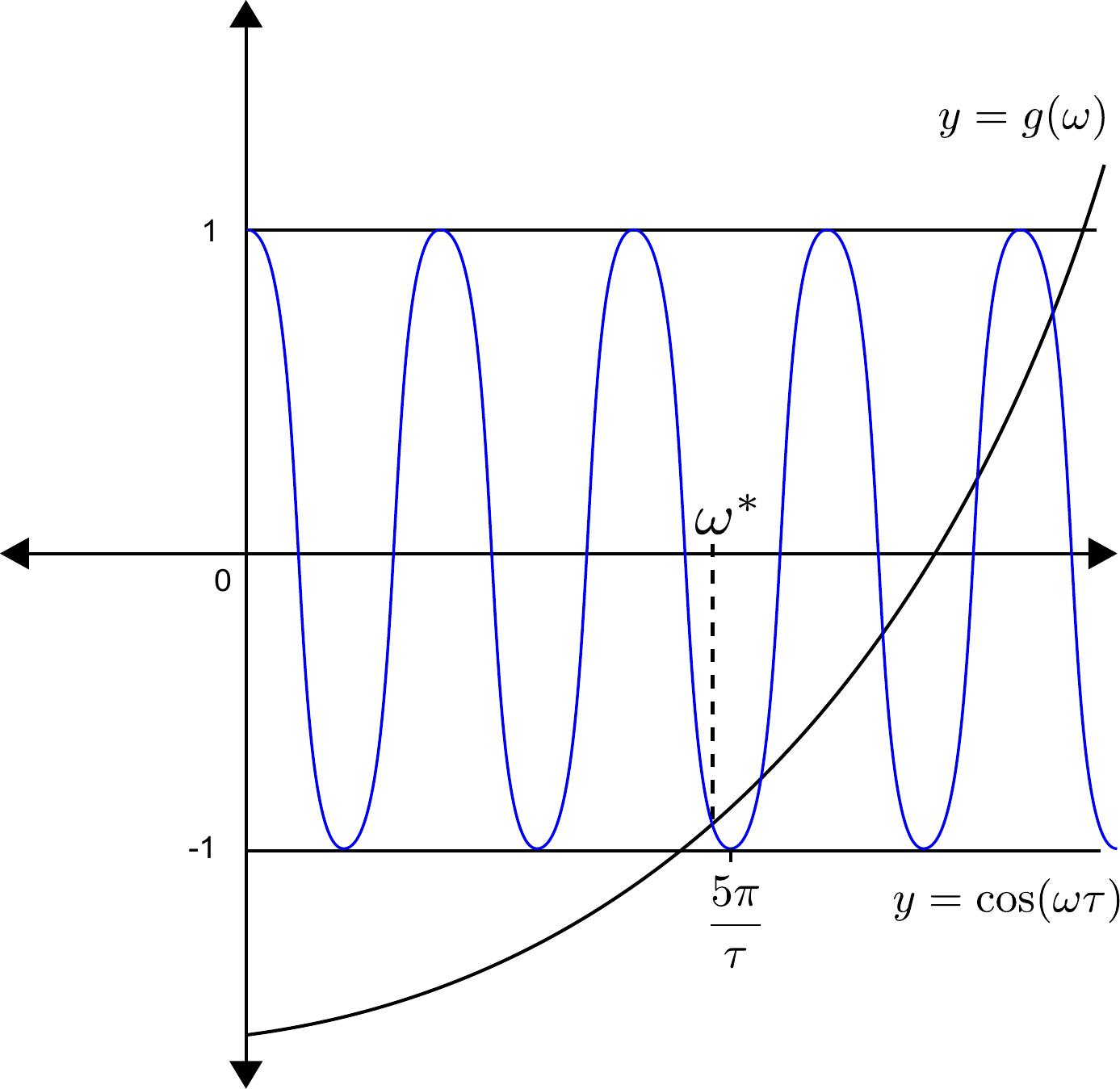}}}
      \subfigure[$bc<0$  and $|b+c|<|a|<|b-c|$]{\scalebox{0.43}{\includegraphics{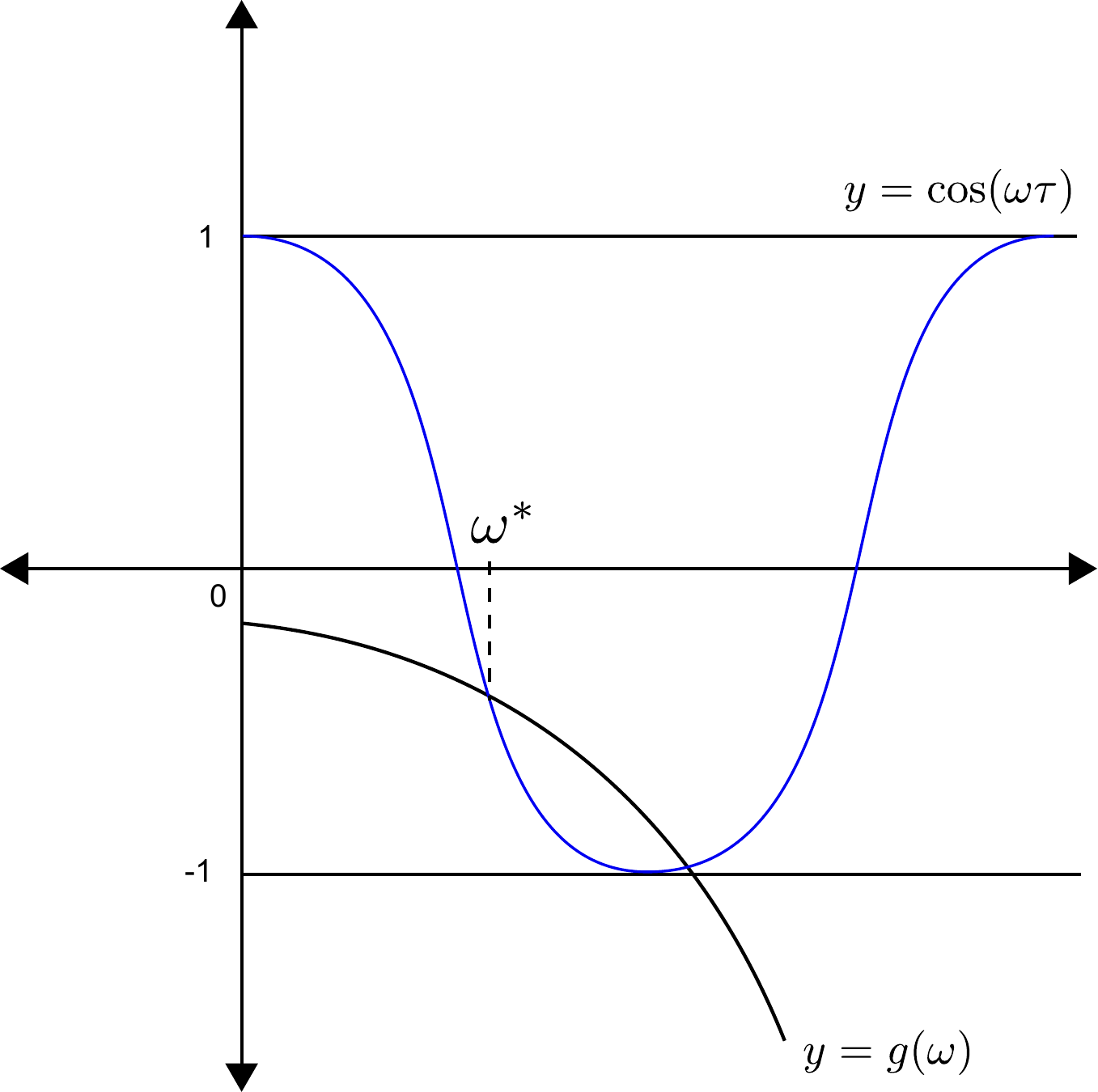}}}
      \subfigure[$bc<0$ and $|a|<|b+c|$]{\scalebox{0.43}{\includegraphics{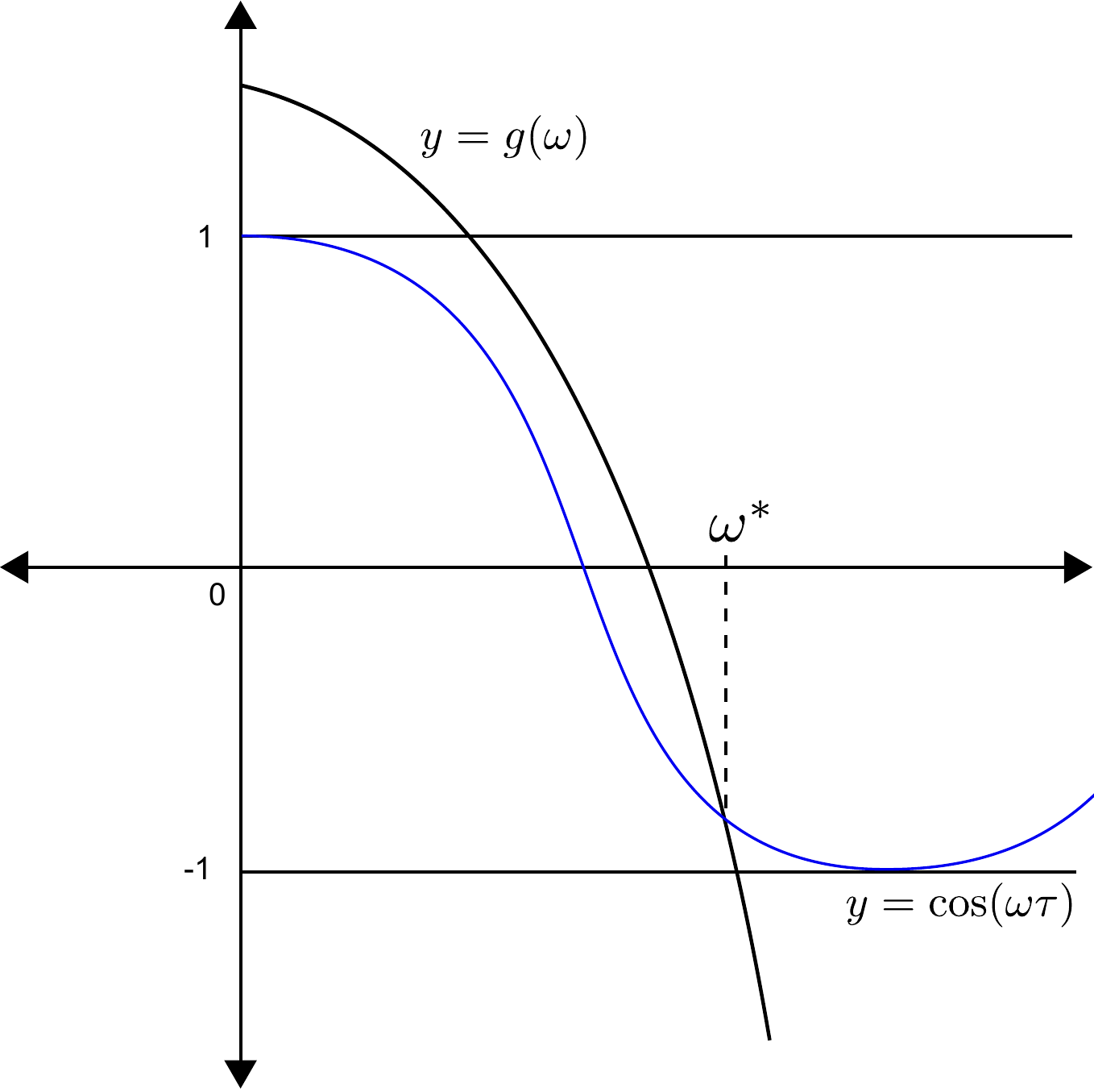}}}
    \subfigure[$bc>0$ and $|a|>|b+c|$]{\scalebox{0.43}{\includegraphics{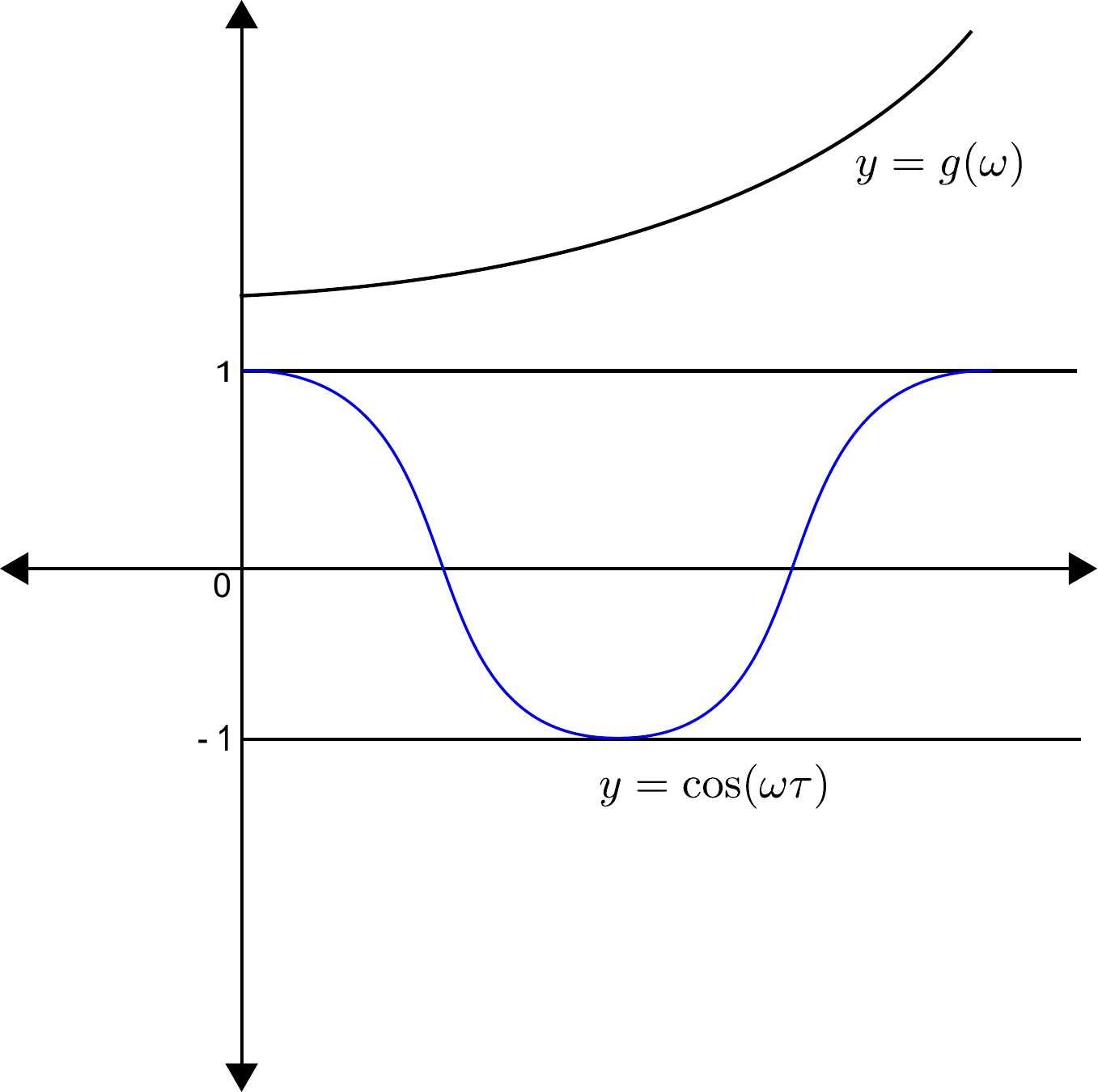}}}
    \subfigure[$bc<0$ and $|a|\geq |b-c|$]{\scalebox{0.43}{\includegraphics{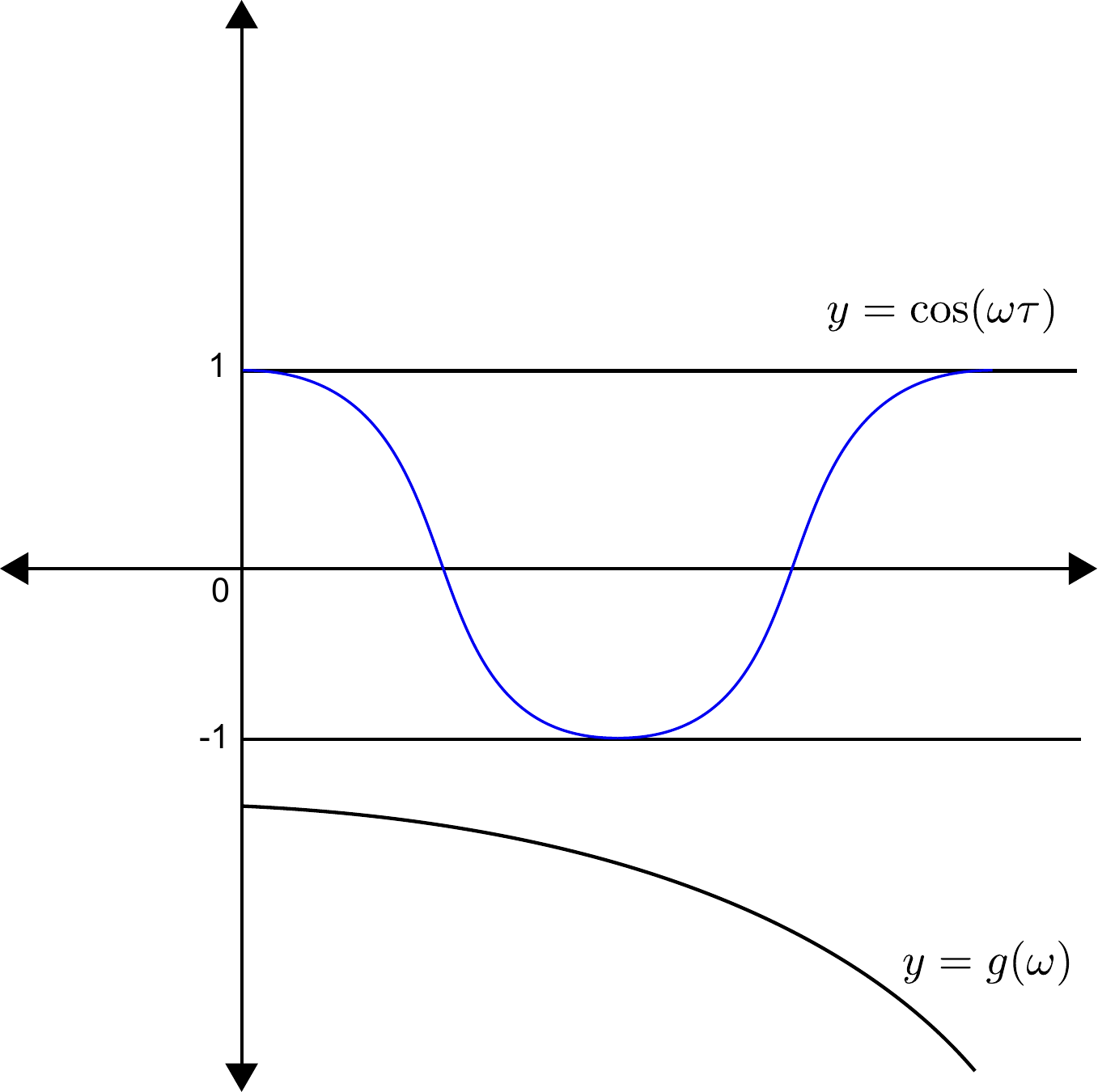}}}
 \caption{~ A first solution $\omega^*$ of \eqref{sistemacaracteristica} as an intersection of the graphs $y=g(\omega)$ and $y=\cos( \omega \tau)$}
 \label{Fig1}
 }
\end{figure}

In Figure \ref{Fig1}, we show the different possibilities given in Lemma \ref{lematau} in terms of the intersections between the curves $y=g(\omega)$ and $y=\cos(\omega\tau)$. 

\begin{remark}
If $\tau\in \mathbb{R}$ is large enough, there could be more than one purely imaginary root of equation \eqref{caracteristica}. For this reason, since we have to study the existence of a pair of simple characteristic roots of \eqref{caracteristica} when \eqref{caracteristica} has no other roots with zero real parts (see, e.g., \cite[Chapt.  8, p.150]{Hale2}), we will slightly strengthen the assumptions of Lemma \ref{lematau} to restrict ourselves to the case where  the uniqueness of the simple root is guaranteed for some pair of delays $(r_0,\sigma_0)$.
\end{remark}

Let $\omega^*$ be the first solution of equation \eqref{equiv} and define 
\begin{equation}\label{taucritico}
\tau^*  \stackrel{def}{=} \min \left\lbrace \tau>0 : \textit{equations \eqref{equiv} and \eqref{tangente} are fulfilled for} \  \omega=\omega^* \right\rbrace, 
\end{equation}
to have the following lemma. 

\begin{lema}\label{lemma23}
Let $\tau\in \mathbb{R}$ and $|\tau|<\tau^*$ be fulfilled, then there exists a pair of delays $(r, \sigma)=(r_0,\sigma_0)$ on the curve $\sigma=r-\tau$ in the $(r, \sigma)$--plane within $\mathbb{R}_+^2$ such that the following assertions are valid: 
\begin{enumerate}
\item[(i)] If $bc>0$ and $ |b-c|\leq |a|<|b+c|$, then there exists a unique pair of simple purely imaginary roots $\pm i\omega^*$ of \eqref{caracteristica} at $(r_0, \sigma_0)$ where $\omega^* \in \left(0, \sqrt{(b+c)^2-a^2}\right]$.
\item[(ii)] If $bc>0$ and $ |a|<|b-c|$, then there exists a unique pair of simple purely imaginary roots $\pm i\omega^*$ of \eqref{caracteristica} at $(r_0, \sigma_0)$ where $\omega^*\in \left[\sqrt{(b-c)^2-a^2}, \sqrt{(b+c)^2-a^2} \right]$.
\item[(iii)] If $bc<0$ and $ |a|<|b+c|$, then there exists a unique pair of simple purely imaginary roots $\pm i\omega^*$ of \eqref{caracteristica} at $(r_0,\sigma_0)$ where $\omega^*\in \left[ \sqrt{(b+c)^2-a^2}, \sqrt{(b-c)^2-a^2} \right]$. 
\end{enumerate}
\end{lema}

\begin{proof}
Assuming $bc>0$ and $|b-c|\leq |a|<|b+c|$,  then $\omega^* \in [0, \omega_1]$ where $g(\omega_1)=1$ and therefore $g(\omega^*)$ take the values in $[g(0),1]$. Put
\begin{equation}\label{pcero}
p(\omega^*) \stackrel{def}{=} g(\omega^*)-\cos(\omega^* \tau).
\end{equation}
Since $p(0)<0\leq p(\omega_1)$, then there exists a solution $\omega^*\in (0, \omega_1]$ to \eqref{pcero}. On the other hand suppose that there exist two solutions to \eqref{pcero}, then $p'(\omega^*)=0$ for some $\omega^*\in (0, \omega_1 ]$ which contradicts the fact that $\omega^*/bc+\tau\sin(\omega^*\tau)>0$ for all $\omega^* \in (0,\omega_1]$ and $|\tau|<\tau^*$, hence the uniqueness of the solution holds and the conclusion follows from Lemma \ref{equivalentes}. The proofs in the other cases are analogous, therefore, they are omitted.
\end{proof}
For $\tau\in \mathbb{R}$ fixed, now we prove the existence of a unique pair of solutions $\lambda=\pm i \omega^*$ to  \eqref{caracteristica} crossing transversally the imaginary axis. To check the \emph{Transversality Condition} we have the following lemma.
\begin{lema}\label{lemmatrans} Assume that the parameters $a,b,c$ and $\tau$ satisfy any of the conditions in Lemma \ref{lemma23} and let $(r_0,\sigma_0,\omega^*)$ denote the values given in that Lemma for some fixed $\tau$. 
Let $\lambda(r)=\mu(r)+i \omega(r)$ be a root of  \eqref{caracteristica} such that $\mu(r_0)=0$ and $\omega(r_0)=\omega^*$. Then $\mu'(r_0)>0$.\end{lema}  

\begin{proof}
 Let $\lambda \in \mathbb{C}$ be a solution of \eqref{caracteristica}, then
\begin{equation*}
\frac{dh(\lambda )}{d\lambda}= 1-b r e^{-\lambda  r}-c \sigma e^{-\lambda \sigma}.
\end{equation*}
For a fixed value of $\tau$, we have 
\begin{equation*}
\frac{dh(i\omega^*)}{d\lambda}= 1+a (r_0-\tau) - b\tau\cos(\omega^* r_0) + i[\omega^*  (r_0-\tau)+ b\tau \sin(\omega^* r_0)].
\end{equation*}
Therefore,
\begin{eqnarray*}
\frac{d}{d\lambda} Re (h(i \omega^*)) &=& 1+a(r_0-\tau) - b\tau \cos(\omega^* r_0), \\
\frac{d}{d\lambda} Im (h(i\omega ^*)) &=& \omega^* (r_0-\tau)  + b\tau \sin(\omega^* r_0).
\end{eqnarray*}
In addition $\lambda (r)$ satisfies 
\begin{equation*}
\lambda(r)+a + e^{-\lambda(r) r} f(\lambda(r)) =0, 
\end{equation*}
where $f(\lambda(r))=b+ce^{\lambda(r) \tau}$. Taking the derivative with respect to $r$ it follows that
\begin{equation*}
\lambda '(r) (1- b\tau e^{-\lambda(r) r}+(\tau-r)f(\lambda(r))e^{-\lambda(r) r}) = \lambda(r) f(\lambda (r)) e^{-\lambda (r)r}.
\end{equation*}
From \eqref{caracteristica} we have that $f(\lambda(r))e^{-\lambda(r) r}=-(\lambda(r)+a)$ and, therefore 
\begin{equation*}
\lambda '(r)= - \frac{\lambda (r)(\lambda(r)+a)}{1-b\tau e^{-\lambda(r)r}+ (r-\tau)(\lambda(r)+a)}.
\end{equation*}
Since $\lambda(r)=\mu (r)+i\omega (r)$ is a root of \eqref{caracteristica} such that $\mu(r_0)=0$ and $\omega (r_0)= \omega^*$, then
\begin{equation*}
\mu'(r_0) = \frac{\omega^*(\omega^*-b\tau (\omega^*\cos(\omega^*r_0)+ a\sin(\omega^*r_0)))}{(1+a (r_0-\tau) - b\tau \cos(\omega^* r_0))^2+( \omega^*  (r_0-\tau) + b\tau \sin(\omega^* r_0))^2}.
\end{equation*}
From \eqref{sistemaimaginariapura} it follows that 
\begin{equation*}
\mu'(r_0)= \frac{\omega^*(\omega^*+bc \  \tau  \sin( \omega^* \tau))}{(1+a  (r_0-\tau) - b\tau \cos(\omega^* r_0))^2+( \omega^* (r_0-\tau)+ b\tau\sin(\omega^* r_0))^2}.
\end{equation*}
Under any of the assumptions of Lemma \ref{lemma23} we have that $\omega^*+bc \  \tau  \sin( \omega^* \tau)>0$, which leads us to conclude that $\mu'(r_0)>0$. 
\end{proof}

In order to study the Hopf bifurcation of \eqref{ecfinal} by choosing one of the delays as the bifurcation parameter, specifically $r$, we need to avoid the threesomes $(r_0, \sigma_0,\omega^*)$ given in Lemma \ref{lematau} that satisfy $r_0\sigma_0=0$. To this end (see, e.g., \cite[ Thm. 4.7 (b), p.53]{Hsmith}) we consider the following assumptions 
\begin{description}
\item[(C.1)] $a>c-b$,
\item[(C.2)] $a>b-c$. 
\end{description}
Note that if conditions  \textbf{(A)},  \textbf{(B)} and \textbf{(C.1)} or  \textbf{(C.2)} hold, then the number of right half plane zeros of \eqref{caracteristica} is zero for $r=\sigma=0$ and with the delays $(r,\sigma) \in \mathbb{R}^2_+$ varying in the semi axis $(r,0)$ or the semi axis $(0,\sigma)$ on the $(r, \sigma)$--plane.

Lemma \ref{lemma23} introduces the existence of a unique pair of simple and purely imaginary roots to \eqref{sistemaimaginariapura} for some pair of delays $(r_0,\sigma_0)\in \mathbb{R}^2_+$. The natural question arises: How to calculate the values of the pair $(r_0,\sigma_0)$? To answer this question, note that from Lemma \ref{lematau} it follows that for $\omega^*>0$ there exists a sequence $\lbrace r_k \rbrace_{k\in \mathbb{N}}$ such that $\omega^*$ satisfies \eqref{sistemaimaginariapura}, namely 
 \begin{equation}\label{seconddelay}
 r_k= \frac{1}{\omega^*} \left[ \arccos\left(-\frac{c\omega^*\sin(\omega^*\tau)+a(b+c\cos(\omega^*\tau))}{{\omega^*} ^2+a^2}\right)+2k\pi\right] \qquad (k=0, 1, 2,  \ldots ).
 \end{equation}
Once we have determined $\omega^*$ for $\tau\in \mathbb{R}$ fixed, then formulae \eqref{uvectores} and \eqref{seconddelay} allows us to calculate the pair of delays $(r_0, \sigma_0)$ in the following manner. Without loss of generality assume that the delay $r$ is the bifurcation parameter and define
\begin{equation}\label{rcero}
r_0  \stackrel{def}{=} \min  \lbrace r_k: k\geq k_{\tau} \  \textit{and}  \ \  i \omega^*  \   \textit{is a simple and unique root of \eqref{caracteristica}} \rbrace , \ \ \sigma_0=\overline{\sigma},
\end{equation}
where $r_k, \ \overline{\sigma}$ and $k_{\tau}$ are given in the proof of Lemma \ref{lematau}, i.e., $k_{\tau}$ is such that $r_k=\tau+\sigma_k\geq 0$ for all $k\geq k_{\tau}$, $\bar{\sigma}\in [0,2\pi/\omega^*)$ is the first value such that \eqref{uvectores} holds, and $r_k$ satisfy \eqref{seconddelay}. 

Assume that the parameters $a, b, c$ and $\tau$ satisfy conditions \textbf{(C.1)} or \textbf{(C.2)} together with any of the hypotheses of Lemma \ref{lemma23}. Then, we propose the following simple procedure to find a unique pair of simple characteristic roots $\mu(r)\pm i\omega(r)$ to \eqref{caracteristica} crossing transversally the imaginary axis at a first point $r=r_0$.

For $\tau\geq 0$:
\begin{enumerate}
\item[(1)] Form the curve $\sigma=r-\tau$ in the $(r, \sigma)$--plane within $\mathbb{R}_+^2$ and initiating from the point $( \tau,0)$ find the first values $(r_0,\sigma_0)$ given in \eqref{rcero}. 
\item[(2)]In Lemma \ref{lemma23} take $(r,\sigma)=(r_0, \sigma_0)$ which are such that \eqref{caracteristica} has a first zero crossing transversally the imaginary axis and $\sigma_0=r_0-\tau$ . 
\end{enumerate}

Analogously for $\tau<0$.

\begin{remark}
If the condition \textbf{(C.1)} is fulfilled we can apply the above procedure for $\tau>0$, on the other hand if  \textbf{(C.2)} is fulfilled we can apply the procedure for $\tau<0$. Although the main importance of our results is that they are applicable in the case $\tau\neq 0$, note that the assumptions of our main results do not exclude the case $\tau=0$. Indeed, for $\tau=0$ the equation  \eqref{caracteristica} is
\begin{equation}\label{caracttau0}
\lambda + a + (b+c)e^{-\lambda r}=0.
\end{equation}
If $\lambda=i\omega^*$, then  \eqref{sistemaimaginariapura} becomes
\begin{gather*}
\begin{aligned}
a &=-(b+c)\cos (\omega^* r), \\
\omega^* &= \ \ (b+c) \sin (\omega^* r),
\end{aligned}
\end{gather*}
which implies that $\omega^*=\pm \sqrt{(b+c)^2-a^2}$. If $|a|<|b+c|$ holds, then \eqref{caracttau0} has a pair of imaginary roots $\pm i\omega^*$ at the sequence 
\begin{equation*}
r_k=\frac{1}{\sqrt{(b+c)^2-a^2}}\left [ \arccos \left( -\frac{a}{b+c} \right) +2k\pi \right]  \qquad (k=0, 1, 2,  \ldots ).
\end{equation*}
If $r=0$ we have that $\lambda=-(a+b+c)$, therefore, if conditions $\textbf{(A)}$ and $|a|<|b+c|$ hold, then the equation \eqref{ecfinal} undergoes a Hopf bifurcation at the zero equilibrium when $r=\sigma=r_0$, i.e., if $r=\sigma$ then the zero equilibrium of equation \eqref{ecfinal} is asymptotically stable for $0<r<r_0$ and unstable for $r>r_0$. The above analysis corresponds to already known results (see, e.g., \cite{Gopalsamy3,Hsmith}). 
\end{remark}

According to Lemmas \ref{lemma23}, \ref{lemmatrans} and by using Rouche's Theorem \cite{Dieudonne}, now on account of \eqref{taucritico} and \eqref{rcero}, we formulate our main result dealing with the stability and local Hopf bifurcations for the zero equilibrium of the equation \eqref{ecfinal}.

\begin{teo}\label{stability}
The following assertions are valid:
\begin{enumerate}
\item[(i)] If the parameters $a,b$ and $c$ satisfy \textbf{(A)}, $\ bc>0$ and $ \ |a|> |b+c|$, then the zero equilibrium of equation \eqref{ecfinal} is locally asymptotically stable for any delays $r\geq 0$ and $\sigma \geq 0$.
\item[(ii)] If the parameters $a,b$ and $c$ satisfy \textbf{(A)}, $\ bc<0$ and $ \ |a|\geq |b-c|$, then the zero equilibrium of equation \eqref{ecfinal} is locally asymptotically stable for any delays $r\geq 0$ and $\sigma \geq 0$.
\item[(iii)] Assume that the parameters $a,b$ and $c$ satisfy \textbf{(A)}, $\ bc<0$ and $\ |b+c|<|a|<|b-c|$. Let, in addition, $\tau\in \mathbb{R}$ be such that $|\tau|<\tau^*$. Then, the zero equilibrium of equation \eqref{ecfinal} is locally  asymptotically stable for any delays $r\geq 0$ and $\sigma \geq 0$ such that $\tau=r-\sigma$.  
\end{enumerate}
\end{teo}

\begin{teo}\label{bifurcation}
There exists a bifurcation parameter $r_0>0$ such that the following assertions are valid: 
\begin{enumerate}
\item[(i)] Assume that the parameters $a,b$ and $c$ satisfy  \textbf{(A)}, $bc>0, \ |b-c|\leq |a| < |b+c|$, and \textbf{(C.1)} or \textbf{(C.2)}. Let, in addition, $\tau\in \mathbb{R}$ be such that $|\tau|<\tau ^ *$. Then, for $r\in[0, r_0)$ the trivial solution to \eqref{ecfinal} is asymptotically stable and the equation \eqref{ecfinal} undergoes a Hopf bifurcation at the trivial solution when $r=r_0$. 
\item[(ii)] Assume that the parameters $a,b$ and $c$ satisfy \textbf{(A)}, $bc>0, \ |a|<|b-c|$, and \textbf{(C.1)} or \textbf{(C.2)}. Let, in addition, $\tau\in \mathbb{R}$ be such that $|\tau|<\tau ^ *$. Then, for $r\in[0, r_0)$ the trivial solution to \eqref{ecfinal} is asymptotically stable and the equation \eqref{ecfinal} undergoes a Hopf bifurcation at the trivial solution when $r=r_0$. 
\item[(iii)] Assume that the parameters $a,b$ and $c$ satisfy \textbf{(A)}, $ bc<0, \ |a|<|b+c|$, and \textbf{(C.1)} or \textbf{(C.2)}. Let, in addition, $\tau\in \mathbb{R}$ be such that $|\tau|< \tau^ *$.Then, for $r\in[0, r_0)$ the trivial solution to \eqref{ecfinal} is asymptotically stable and the equation \eqref{ecfinal} undergoes a Hopf bifurcation at the trivial solution when $r=r_0$. 
\end{enumerate}
\end{teo} 

\section{Direction and stability of the Hopf bifurcation}\label{S3}
In this section, by applying normal form theory, we study the direction of the Hopf bifurcation and the stability of the bifurcating periodic solutions. We refer to \cite{Faria, Hale2} for the results and notations involved. 

Throughout this section, we will assume that the conditions and conclusions of Theorem \ref{bifurcation} are fulfilled, therefore, the delay $\sigma$ in \eqref{ecfinal} will no longer be interpreted as a free parameter in the sense that the difference between the delays is assumed constant. Now we introduce some preliminary notation and sets that we use throughout this section.

Let $\kappa>0$ and define the phase space $C\stackrel{def}{=}C([-\kappa , 0], \mathbb{C})$ equipped with the sup norm. Consider the following Banach space
\begin{equation*}
BC=\lbrace \phi:[-\kappa, 0] \rightarrow \mathbb{R}: \phi \ \textit{is continuous on} \ [-\kappa, 0), \exists \lim_{\theta \rightarrow 0^-} \phi(\theta) \in \mathbb{R} \rbrace.
\end{equation*}
Without loss of generality, let $\max(r,r-\tau)=r$ and define the linear functional $L$ on $BC$ as
\begin{equation*}
L(r)\phi = \int_{-r}^0 \phi (\theta) d\eta(\theta)=-a\phi(0)-b\phi(-r)-c\phi(\tau-r),
\end{equation*}
where 
\begin{eqnarray*}
\eta (\theta) = \left\lbrace
\begin{array}{lcc}
0 & \mbox{if} & \ \qquad \qquad \theta = -r, \\
-b & \mbox{if} &\ \ \qquad  -r < \theta \leq \tau-r, \\
-(b+c) &\mbox{if} & \tau-r <\theta \leq 0, \\
-(a+b+c) & \mbox{if} &  \ \ \ \ \ \qquad \theta >  0.
\end{array}
\right.
\end{eqnarray*}
Introducing the new parameter $\alpha=r-r_0$, we will consider $$L(\alpha)\phi= -a\phi(0)-b\phi(-(r_0+\alpha))-c\phi( \tau -(r_0+\alpha)).$$ 
Due to the inclusion of $\mathbb{R}^n$ into $C$, instead of \eqref{ecfinal} we consider the equation 
\begin{equation*}
\dot{u}(t)=L(\alpha)u_t+F(u_t, \alpha),
\end{equation*}
where $\alpha\in V,~ V$ a neighborhood of zero in $\mathbb{R}$ and $F(u_t, \alpha)=f(u_t, \alpha)-L(\alpha)u_t$ with $f(u_t, \alpha)=f(u(t), u(t-(r_0+\alpha)), u(t-(r_0+\alpha-\tau)))$. In order to calculate the normal forms up to the third order, we will assume that the functions $L: V \rightarrow \mathcal{L}(C; \mathbb{R})$ and $f: C \times V \rightarrow \mathbb{R}$ are sufficiently differentiable such that $F: C\times V \rightarrow \mathbb{R}$ is a $C^3$ function. 

Let $L_0=L(0), \ \Lambda = \lbrace  i\omega^*,- i \omega^* \rbrace$ and $P$ be the center space of $y'(t)=L_0y_t$. Decomposing $C$ by $\Lambda$ as $C=P\oplus Q$, we choose bases $\Phi, \Psi$ for $P$ and $P^*$ respectively as
\begin{eqnarray*}
P = \textit{span} \Phi,  && \Phi(\theta)=(\varphi_1(\theta), \varphi_2(\theta))= (e^{i\omega^* \theta}, e^{-i\omega^*\theta}), \ \ - r \leq \theta \leq 0,  \\
 P^*=\textit{span} \Psi, && \Psi(s) = \textit{col} (\psi_1(s), \psi_2(s))=\textit{col} (\psi_1(0)e^{-i\omega^*s}, \overline{\psi_1(0)}e^{i\omega^* s}), \\ && \qquad \qquad 0\leq s \leq r,
\end{eqnarray*}
where 
\begin{equation*}
\psi_1(0)\stackrel{def}{=} \left(1-L_0\left(\theta e^{i\omega^*\theta}\right) \right)^{-1}=\left(1-b \tau e^{-i\omega^*r_0}+(r_0-\tau)(i\omega^*+a)\right)^{-1}.
\end{equation*}
Writing the Taylor expansions $L(\alpha)=L_0+L_1(\alpha)+\frac{1}{2}L_2(\alpha)+ h.o.t. , \ F(v, \alpha)=\frac{1}{2}F_2(v, \alpha)+ \frac{1}{3!} F_3(v, \alpha)+h.o.t.,$  where h.o.t means higher order terms and  $L_j, F_j$ are the $j$th Fr\'echet derivative of $L$ and $F$ in the variables $\alpha$ and $(v,\alpha)$, respectively, we will consider the Taylor formulas 

\begin{equation*}
\begin{array}{l} 
F_2( \phi, 0) = a_{11}\phi^2(0)+a_{22} \phi ^2(-r_0) +a_{33}\phi^2(-(r_0-\tau )) \\
\qquad  \qquad  \ \   +2a_{12}\phi(0)\phi(-r_0)+ 2a_{13}\phi(0)\phi(-(r_0-\tau ))\\
\qquad  \qquad \ \   + 2a_{23}\phi(-r_0)\phi(-(r_0-\tau )), \\

\\

 F_3( \phi, 0) = b_{111}\phi^3(0)+b_{222}\phi^3(-r_0)+b_{333}\phi^3(-(r_0-\tau )) \\
 \qquad \qquad   \ \   + 3b_{112}\phi^2(0)\phi(-r_0)+3b_{113}\phi^2(0)\phi(-(r_0-\tau )) \\ 
\qquad  \qquad \ \ +3b_{122}\phi(0)\phi^2(-r_0)+3b_{133} \phi(0)\phi^2(-(r_0-\tau )) \\
 \qquad \qquad \ \   + 6b_{123}\phi(0)\phi(-r_0)\phi(-(r_0-\tau )) \\ 
\qquad  \qquad \ \  + 3b_{223}\phi^2(-r_0)\phi(-(r_0-\tau )) \\
\qquad \qquad \ \   + 3b_{233}\phi(-r_0)\phi^2(-(r_0-\tau )).
\end{array}
\end{equation*}
Let 
\begin{eqnarray*}
E_1 &=& 3 \psi_1(0) \cdot \Bigl( b_{111}+b_{222}e^{-i\omega^*r_0} +b_{333}e^{i\omega^*(\tau-r_0)}+b_{112}(e^{i\omega^*r_0}+2e^{-i\omega^*r_0})\\  && + b_{113}(e^{-i\omega^*(\tau-r_0)}+2e^{-i\omega^*(r_0-\tau)}) +b_{122}(2+e^{-2i\omega^*r_0}) \\ &&+b_{133} (2+e^{-2i\omega^*(r_0-\tau)}) +b_{123}(2e^{-i\omega^*\tau}+2e^{i\omega^*\tau}+2e^{i\omega^*(\tau-2r_0)}) \\
&&+ b_{223}(e^{-i\omega^*(r_0+\tau)}+2e^{i\omega^*(\tau-r_0)})+b_{233}(2e^{-i\omega^*r_0}+2e^{i\omega^*(2\tau-r_0)}\Bigr), \\
E_2 &=&\Bigl(a_{11}+a_{22}+a_{33}+2a_{12}\cos(\omega^*r_0)+2a_{13}\cos(\omega^*(r_0-\tau)) \\
&& +2a_{23}\cos(\omega^*\tau)\Bigr)\cdot \Bigl(a+b+c\Bigr)^{-1}, \\
E_3 &=& \psi_1(0) \cdot \Bigl( a_{11}+a_{12}+a_{13}+(a_{12}+a_{22}+a_{23})e^{-i\omega^*r_0} \\
&& +(a_{13}+a_{23}+a_{33})e^{i\omega^*(\tau-r_0)}\Bigr) , \\
E_4 &=& \psi_1(0)\cdot \Bigl( a_{11}+a_{22}e^{-2i\omega^*r_0} +a_{33}e^{2i\omega^*(\tau-r_0)}+ 2a_{12}e^{-i\omega^*r_0} \\
&& +2a_{13}e^{i\omega^*(\tau-r_0)}+2a_{23}e^{i\omega^*(\tau-2r_0)}\Bigr)\cdot \Bigl( a_{11}e^{2i\omega^*r_0} +a_{22}e^{i\omega^*r_0}\\
&&+a_{33} e^{i\omega^*(r_0+\tau)}+a_{12}(1+e^{3i\omega^*r_0})+a_{13} (e^{i\omega^*(3r_0-\tau)}+e^{2i\omega^*\tau}) \\
&&+a_{23} (e^{i\omega^*(r_0-\tau)}+e^{i\omega^*(r_0+2\tau)})\Bigr) \cdot \Bigl( (a+2i\omega^*)e^{2i\omega^*r_0}+b+ce^{2i\omega^*\tau} \Bigr) ^{-1}.
\end{eqnarray*}
According to \cite{Faria, Hale2}, we can calculate the normal forms on the center manifold to \eqref{ecfinal} as follows.
\begin{prop}
A normal form of  \eqref{ecfinal} on the center manifold of the origin is given by
\begin{equation}\label{formanormal}
\dot{x}= B x+  \left(
\begin{array}{c}
B_1x_1\alpha  \\
\overline{B_1}x_2\alpha  \\
\end{array}
\right) +\left(
\begin{array}{c}
B_2x_1^2x_2 \\
\overline{B_2}x_1x_2^2 \\
\end{array}
\right)+ O(|x|\alpha ^2+|x|^4),
\end{equation}
where $B_1, B_2\in \mathbb{C}$ and $B=\mbox{diag} (i\omega^*, -i\omega^*)$.
\end{prop}
The change of variables $w$, where $x_1=w_1-iw_2, x_2=w_1+iw_2$, and the use of polar coordinates $ w_1=\rho\cos(\beta), w_2=\rho \sin(\beta)$ transforms \eqref{formanormal} into
\begin{eqnarray}\label{formanormalcordenadaspolares}
\left\lbrace
\begin{array}{lcl}
\dot{\rho} &=& K_1 \alpha  \rho + K_2\rho ^3+O(\alpha ^2\rho + |(\rho,\alpha )|^4), \\
\dot{\beta} &=& - \omega^* + O(|(\rho, \alpha )|),
\end{array}
\right.
\end{eqnarray}
where 
\begin{equation*}
\begin{array}{l}K_1\stackrel{def}{=} Re\left( \psi_1(0) \cdot (\omega^*(\omega^*-ia) \right))=\mu'(r_0)>0 ,\\
K_2\stackrel{def}{=} \dfrac{1}{6}Re(E_1)+ E_2 \cdot Re(E_3)+ \dfrac{1}{2} Re(E_4).
\end{array}
\end{equation*}

By applying the classical Hopf-bifurcation theory \cite{ChowHale}, now we formulate our main result dealing with the direction of the Hopf bifurcation and the stability of the bifurcating periodic solutions.  
\begin{teo}\label{estabilidaddireccion}
Suppose that any of the assumptions of Theorem \ref{bifurcation} are fulfilled, then the dynamics of \eqref{ecfinal} near the origin is described by  \eqref{formanormalcordenadaspolares}. Since $K_1>0$, then a supercritical Hopf bifurcation occurs at $r=r_0$. Moreover if $K_2< 0$ (resp. $K_2>0$) the bifurcating periodic solutions are asymptotically stable (resp. unstable) on the center manifold. Further, the period of the bifurcating periodic solution is determined by $ T(\varepsilon) =  2\pi / \omega^*(\varepsilon)$ for some $|\varepsilon|<\varepsilon_0$.
\end{teo}

\section{Applications to population models}\label{S4}
In this section, we apply the results obtained in Sections \ref{S2} and \ref{S3} to prove the existence of bifurcating periodic solutions to \eqref{Nicholson},  and the equation \eqref{Intro3} with a delayed harvesting term. The results are illustrated by numerical simulations.

\subsection{Nicholson's model}
We consider the Nicholson's blowflies model \eqref{Nicholson} described in Section \ref{S1}. Here
\begin{equation*}
\delta >0, \ H>0, \ P>0, \ r\geq 0, \ \sigma\geq 0,
\end{equation*}
and we consider  \eqref{Nicholson} subject to the following non-negative initial condition and positive initial value 
\begin{equation}\label{initialcondition}
x(t)=\varphi(t), \ \varphi(t)\geq 0, \ -\max\lbrace r,\sigma\rbrace \leq t \leq 0, \ \ x(0)=x_0>0. 
\end{equation}
\begin{remark}
It is well-known that if $\sigma=0$, then every solution $x$ to the initial value problem \eqref{Nicholson}, \eqref{initialcondition} is positive for $t\geq 0$. However, if $\sigma>0$ there exists a non-negative initial condition such that the solution to the initial value problem \eqref{Nicholson}, \eqref{initialcondition}  becomes negative at some $t_0>0$ (see, e.g., \cite{Wei2021}). Therefore, since our context is biological, below we choose the  appropriate initial conditions to present a numerical simulation of two positive solutions to \eqref{Nicholson} (see Figure \ref{Fig2} below). 
\end{remark}

According to Theorem \ref{stability}, now we formulate a stability criterion for the zero equilibrium of the model \eqref{Nicholson}.

\begin{teo} The following assertions are valid: 
\begin{enumerate}
\item[(i)] If $~0<H<\delta$ and $~ P\leq \delta-H$, then the zero equilibrium of equation \eqref{Nicholson} is locally asymptotically stable for any delays $r\geq 0$ and $\sigma\geq 0$. 
\item[(ii)] Assume that $~ 0<H\leq \delta$ and $~ \delta-H<P<\delta+H$, or $~ 0<\delta<H$ and $~ H-\delta<P<\delta+H$. Let, in addition, $\tau\in \mathbb{R}$ be such that $|\tau|<\tau^*$. Then, the zero equilibrium of \eqref{Nicholson} is locally asymptotically stable for any delays $r\geq 0 $ and $\sigma\geq 0$ such that $\tau=r-\sigma$. 
\end{enumerate}
\end{teo}

If $~ P\leq \delta +H$, then \eqref{Nicholson} has only the zero equilibrium. If $~ P>\delta + H$, then the zero equilibrium is unstable and there exists the following positive equilibrium of \eqref{Nicholson}
\begin{equation*}
x^*=\ln \left( \frac{P}{\delta + H} \right).
\end{equation*}
Therefore, the condition $P>\delta +H$ will be assumed. Furthermore, in order to apply the results of the paper, we will assume that the delays in \eqref{Nicholson} are such that $\sigma=r-\tau$ with $\tau\in \mathbb{R}$. 

Let $x(t)=x^*+u(t)$, then \eqref{Nicholson} becomes 
\begin{equation}\label{nicholsonlineal}
\dot{u}(t)=-\delta u(t)-Hu(t-(r-\tau))+(\delta + H)[u(t-r)e^{-u(t-r)}+x^*(e^{-u(t-r)}-1)]. 
\end{equation}
The linearization of \eqref{nicholsonlineal} around $u=0$ is
\begin{equation*}
\dot{y}(t)=-\delta y(t)-(x^*-1)(\delta + H)y(t-r)-Hy(t-(r-\tau)).
\end{equation*}
According to Theorems \ref{stability}, \ref{bifurcation} and \ref{estabilidaddireccion}, now we  formulate our main result dealing with the stability of the positive equilibrium of \eqref{Nicholson}, and the existence of local Hopf bifurcations at $x=x^*$.

\begin{teo}
The following assertions are valid:
\begin{enumerate}
\item[(i)] If $~ 0<H< \delta$ and $~ 1< x^* < 2\delta / (\delta +H)$, then the positive equilibrium $x^*>0$ of equation \eqref{Nicholson} is locally asymptotically stable for any delays $r\geq 0$ and $\sigma \geq 0$. 
\item[(ii)] If $~ 0<H<\delta$ and $~ 2H/(\delta+H) \leq x^*<1$, then the positive equilibrium $x^*>0$ of equation \eqref{Nicholson} is locally asymptotically stable for any delays $r\geq 0$ and $\sigma \geq 0$. 
\item[(iii)] Assume that $~ 0<H\leq \delta$ and $~ x^*< 2H/(\delta+H)$, or $~ 0<\delta<H$ and $~ x^*< 2\delta/(\delta +H)$. Let, in addition, $\tau\in \mathbb{R}$ be such that $|\tau|<\tau^*$. Then, the positive equilibrium $x^*>0$ of equation \eqref{Nicholson} is locally   asymptotically stable for any delays $r\geq 0$ and $\sigma \geq 0$ such that $\tau=r-\sigma$. 
\end{enumerate}
\end{teo}

\begin{teo}\label{Nicholsonbifsigma}
There exists a bifurcation parameter $r_0>0$ such that the following assertions are valid:
\begin{enumerate}
\item[(i)] Assume that $~ 0<H\leq \delta$ and $~ 2\delta/(\delta+H) <x^*\leq 2$, or $~ 0<\delta<H$ and $ ~ 2H/(\delta+H)\leq x^* \leq 2$. Let, in addition, $\tau\in \mathbb{R}$ be such that $|\tau|<\tau^*$. Then, the equation \eqref{Nicholson} undergoes a supercritical Hopf bifurcation at $x=x^*$ when $r = r_0$.
\item[(ii)] Assume that $~ x^*>2$, or $~ 0<\delta<H$ and $~ 1<x^*< 2H/(\delta+H)$. Let, in addition, $\tau \in \mathbb{R}$ be such that $|\tau|<\tau^*$. Then, the equation  \eqref{Nicholson} undergoes a supercritical Hopf bifurcation at $x=x^*$ when $r=r_0$.
\item[(iii)] Assume that $~ 0<\delta<H$ and $~ 2\delta/(\delta+H)<x^*<1$. Let, in addition, $\tau\in \mathbb{R}$ be such that $|\tau|<\tau^*.$ Then, the equation \eqref{Nicholson} undergoes a supercritical Hopf bifurcation at $x=x^*$ when $r=r_0$.
\end{enumerate}
Furthermore, in any of the above cases, if $K_2<0$ (resp. $K_2>0$), then the bifurcating periodic solutions are asymptotically stable (resp. unstable) on the center manifold. 
\end{teo}
Now we consider a particular case of \eqref{Nicholson}, namely the equation 
\begin{equation}\label{simulacionnicholson}
\dot{x}(t)=-2 x(t)-x(t-(r-\tau))+3e^{2.5}x(t-r)e^{-x(t-r)}.
\end{equation}
Here the positive equilibrium is $x^*=2.5$. As usual, by a solution to  \eqref{simulacionnicholson} we understand a continuously differentiable function $x$ that satisfies the problem  \eqref{simulacionnicholson}, \eqref{initialcondition} for $t\geq 0$. According to Theorem \ref{Nicholsonbifsigma} there exists a parameter $r_0>0$ such that all solutions to \eqref{simulacionnicholson} tend asymptotically to the positive equilibrium for every $r<r_0$, and for $r= r_0$ a Hopf bifurcation occurs at $x^*=2.5$. On the other hand, the stability of the bifurcating periodic solutions is determined by the sign of $K_2$. In particular, for $\tau=0.3782,$ we obtain $\omega^*=4.1533, \  r_0=0.5389$ and  $K_2=-0.3573$, i.e., the bifurcating periodic solution to \eqref{simulacionnicholson} is asymptotically stable on the center manifold. In Figure \ref{Fig2} is showed that for $r=r_0$ a supercritical Hopf bifurcation occurs at $x^*=2.5$.

\begin{figure}
\centering
\subfigure[$r=0.45$]{\scalebox{0.3}{\includegraphics{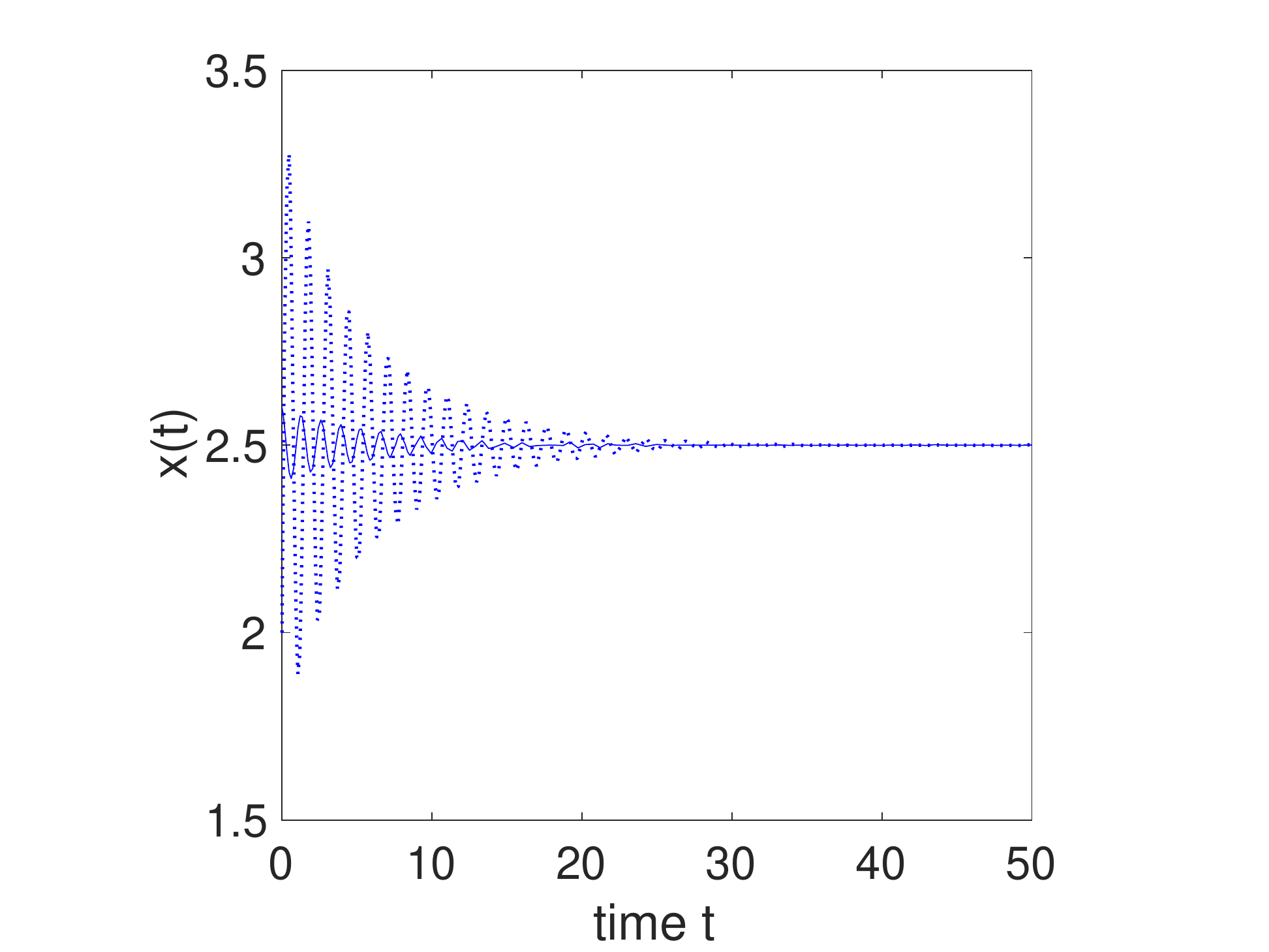}}}
\hspace{0mm}
\subfigure[ $r=0.5$]{\scalebox{0.3}{\includegraphics{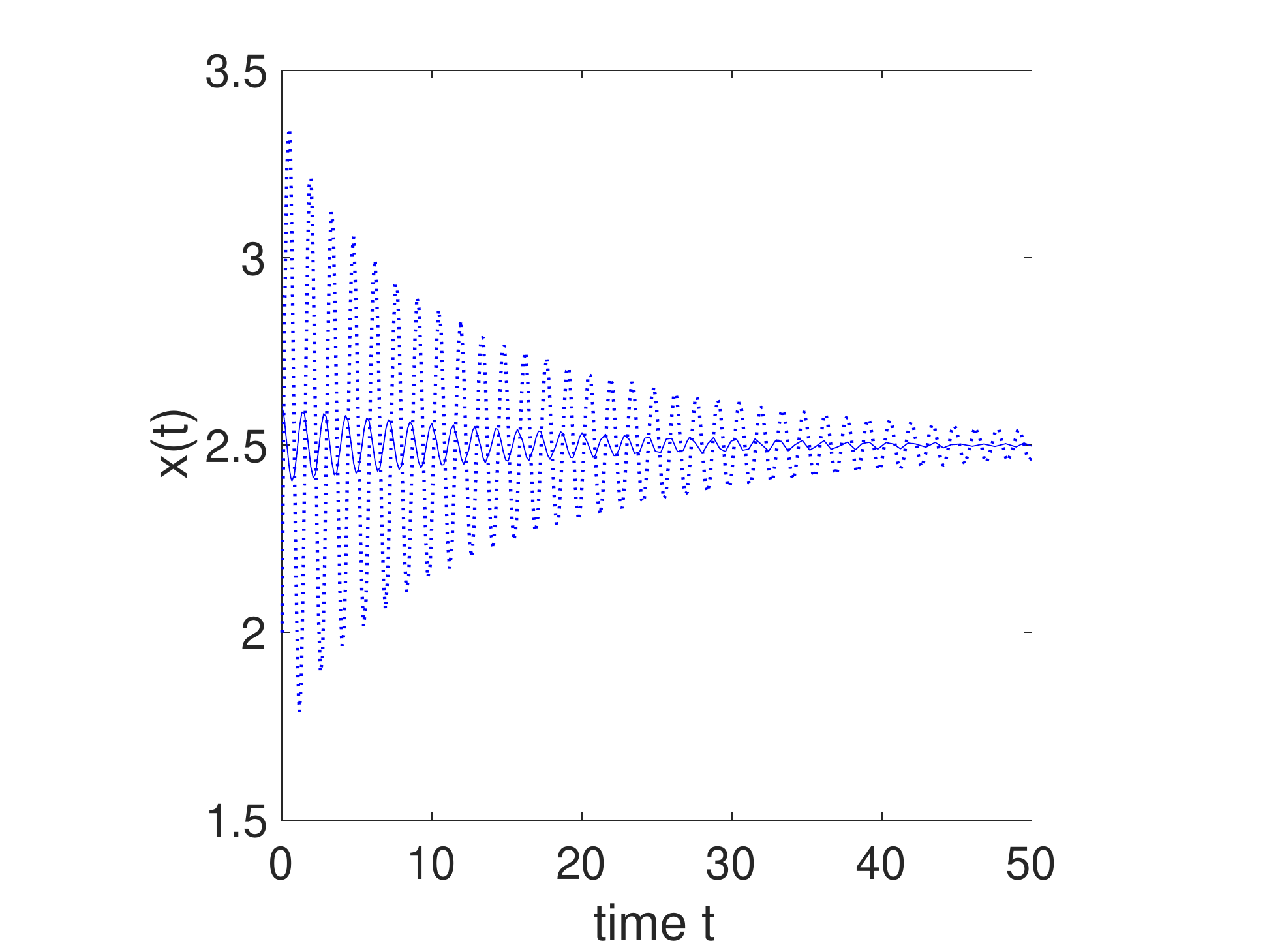}}}
\hspace{0mm}
\subfigure[ $r=r_0$]{\scalebox{0.3}{\includegraphics{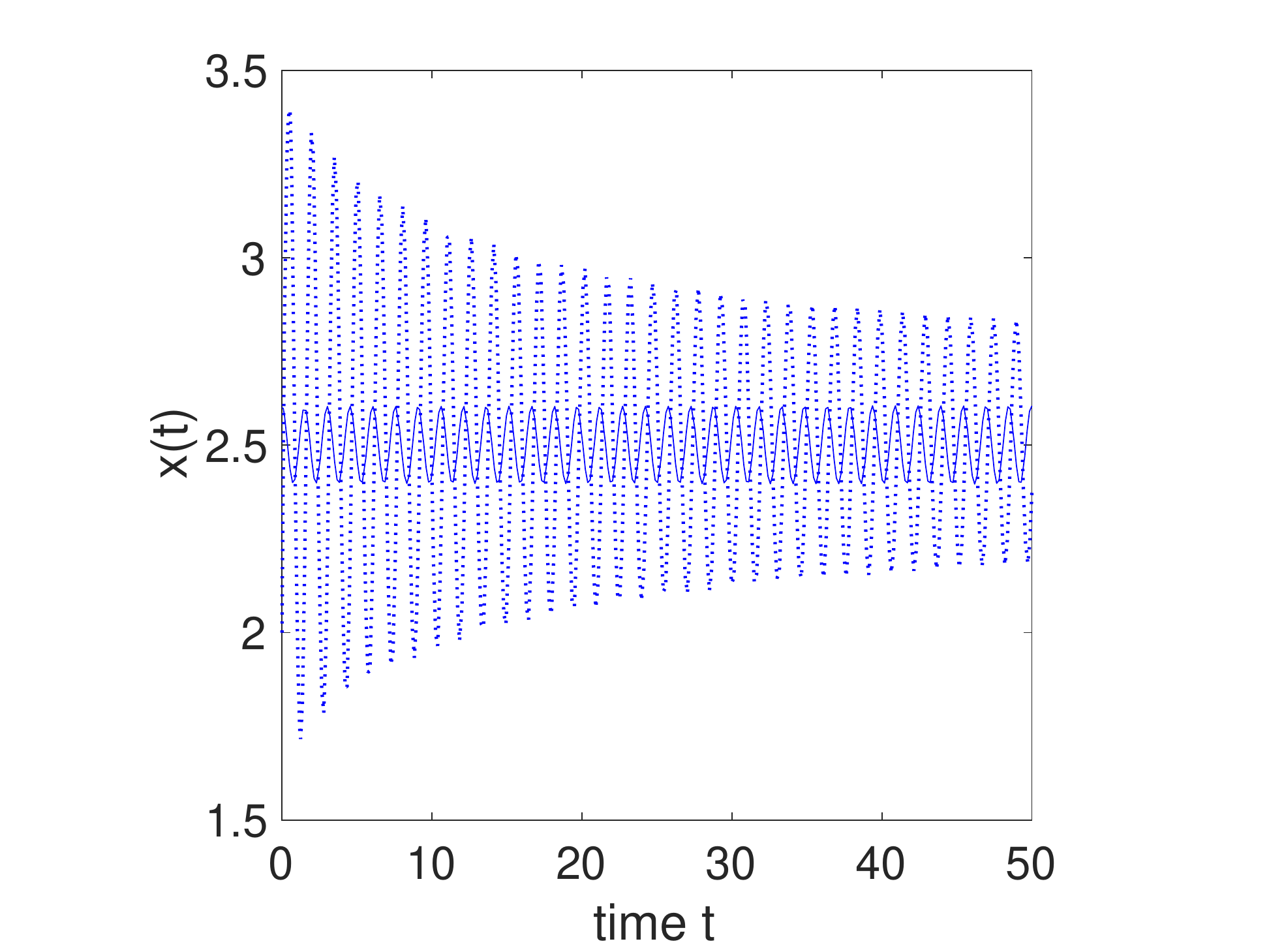}}}
\hspace{0mm}
\subfigure[ $r=0.65$]{\scalebox{0.3}{\includegraphics{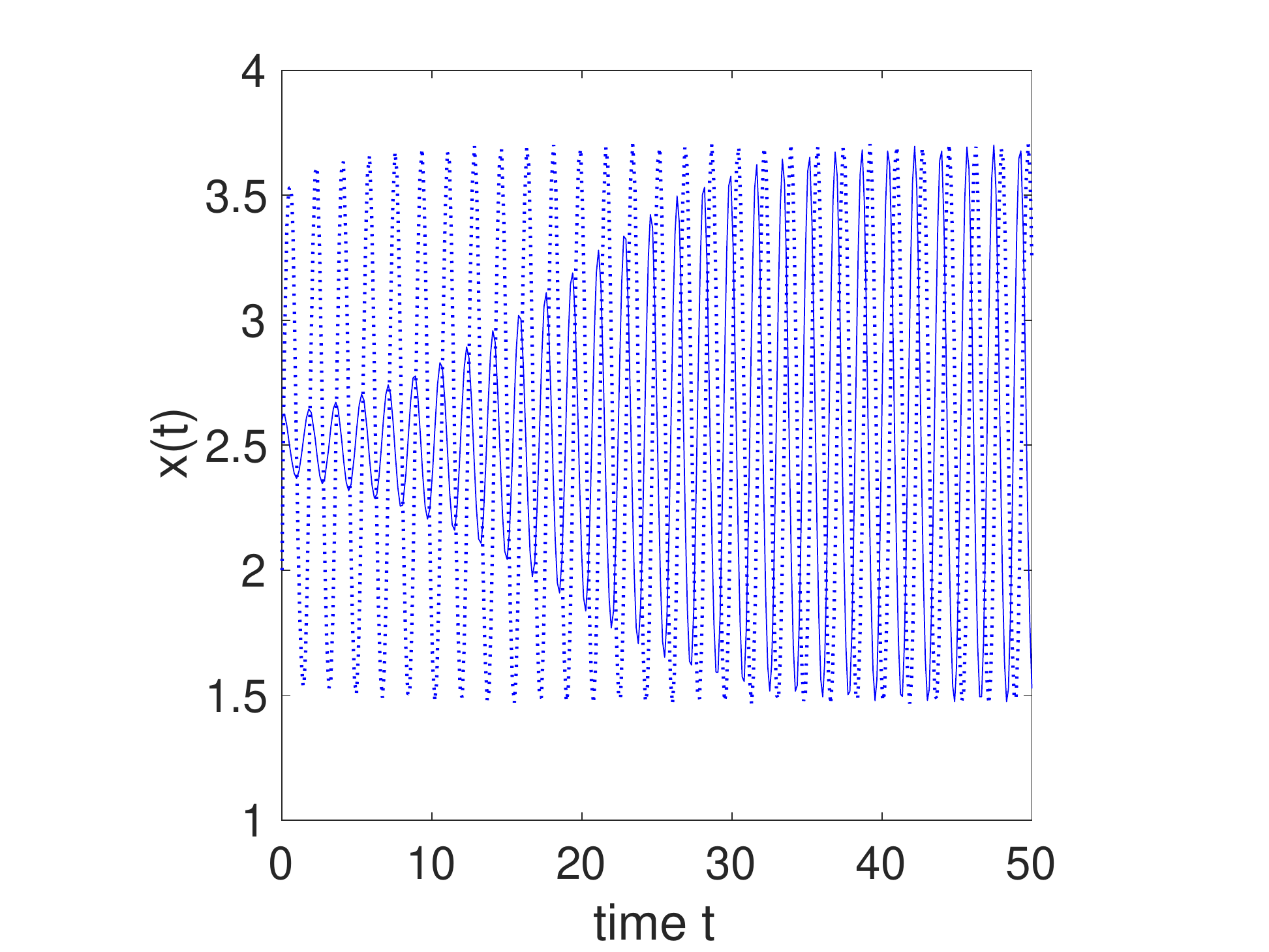}}}
\caption{~ Dynamic behaviour of Nicholson's model \eqref{simulacionnicholson} as $r$ varies, where the initial functions are: $\varphi_1(\theta)=\sin(\theta)+2 $ (dots), $\varphi_2(\theta)=1.3(\cos(\theta)+1)$ (continuous line) for $\theta\in [-r,0]$}
\label{Fig2} 
\end{figure}

\subsection{Mackey-Glass model}

Based on the Mackey-Glass model  \eqref{Intro3} described in Section \ref{S1},  we propose the following Mackey-Glass model with a delayed harvesting term
\begin{equation}\label{MackeyGlass1}
\dot{x}(t)=-\delta x(t)+ \frac{P}{1+x^n(t-r)}-Hx(t-\sigma),
\end{equation}
where 
\begin{equation*}
\delta>0, \ H>0, \ P>0, \ n>0, \ r\geq 0, \ \sigma \geq 0.
\end{equation*}
It is easy to check that \eqref{MackeyGlass1} has a unique positive equilibrium, namely $x=x_*$ which satisfies
\begin{equation*}
x_*^{n+1}+x_*=\frac{P}{\delta+H}.
\end{equation*} 
In \eqref{MackeyGlass1} we assume that $\sigma=r-\tau$ with $\tau\in \mathbb{R}$. Let $x(t)=x_*+u(t)$, then \eqref{MackeyGlass1} becomes
\begin{equation}\label{MackeyGlass2}
\dot{u}(t)=-\delta u(t)-Hu(t-(r-\tau))+\frac{P}{1+(u(t-r)+x_*)^n}-(\delta+H)x_*.
\end{equation}
The linearization of \eqref{MackeyGlass2} around $u=0$ is 
\begin{equation*}
\dot{y}(t)=-\delta y(t)-\frac{Pnx_*^{n-1}}{(x_*^n+1)^2}y(t-r)-Hy(t-(r-\tau)).
\end{equation*}
According to Theorems \ref{stability}, \ref{bifurcation} and \ref{estabilidaddireccion}, we conclude the following.
\begin{teo}\label{MackeyPrincipal}
There exists a bifurcation parameter $r_0>0$ such that the following assertions are valid:
\begin{enumerate}
\item[(i)] If $~Pnx_*^{n-1}/(x_*^n+1)^2 + H <\delta$, then the positive equilibrium $x_*>0$ of equation \eqref{MackeyGlass1} is locally asymptotically stable for any delays $r\geq 0$ and $\sigma \geq 0$.
\item[(ii)] Assume that $~Pnx_*^{n-1}/(x_*^n+1)^2+H>\delta $. Let, in addition, $\tau \in \mathbb{R}$ be such that $|\tau|<\tau^*$. Then, the equation \eqref{MackeyGlass1} undergoes a supercritical Hopf bifurcation at $x=x_*$ when $r=r_0$. Furthermore, if $K_2<0$ (resp. $K_2>0$), then the bifurcating periodic solutions are asymptotically stable (resp. unstable) on the center manifold. 
\end{enumerate}
\end{teo}

According to the notation of Sections \ref{S2} and \ref{S3}, we finish this section presenting in Table \ref{Tab1} the coefficients of the Taylor expansions around the positive equilibrium of the models \eqref{Nicholson} and  \eqref{MackeyGlass1}. These coefficients are useful to determinate the sign of $K_2$ given in Theorems \ref{Nicholsonbifsigma} and \ref{MackeyPrincipal}.

\begin{table}[H]
\centering
\resizebox{9cm}{!} {\begin{tabular}{|p{3cm} |p{3cm} | p{4cm} | p{3cm} |p{1.1cm}|}
 \hline 
Coefficients &  Nicholson's model  & Mackey-Glass model  \\
\hline \hline
$a$ & $\delta$  &  $\delta$  \\
$b$ & $(x^*-1)(\delta+H)$ & $Pnx_*^{n-1}/(x_*^n+1)^2$ \\
$c$ & $H$ & $H$  \\
$a_{11}$  & 0  &0   \\
$a_{22}$ & $(x^*-2)(\delta+H)$& $Pnx_*^{n-2}(1-n+(n+1)x_*^n)/(x_*^n+1)^3$\\
 $a_{33}$ & 0 & 0 \\
$a_{12}$ &0 & 0\\
$a_{13}$ &0 & 0 \\
$a_{23}$ &0 & 0 \\
$b_{111}$ & 0& 0 \\
$b_{222}$ & $(3-x^*)(\delta + H)$ & $Pnx_*^{n-3}((2-n)(n-1)+4(n^2-1)x_*^n-(n+1)(n+2)x_*^{2n})/(x_*^n+1)^4$  \\ 
$b_{333}$ & 0&0 \\
$b_{112}$ & 0&0 \\
$b_{113}$ &0 & 0 \\ 
$b_{122}$ &0 &0 \\ 
$b_{133}$ &0 &0  \\
$b_{123}$ &0 & 0\\
$b_{223}$ &0 & 0\\
 $b_{233}$&0 & 0\\
 \hline 
\end{tabular}}
\caption{ \ \ Coefficients in the Taylor series expansion up  to the third order around the positive equilibrium of Nicholson's model \eqref{Nicholson} and  Mackey-Glass model \eqref{MackeyGlass1}}
\label{Tab1}
\end{table}

\section{Discussion}

In this paper, we study a general class of autonomous scalar differential equations with two discrete delays. We have shown that the difference between the delays plays an important role in determining the stability of the system. At first, we introduce a new method to analyse the distribution of the roots of the corresponding characteristic equation and, by choosing a delay as the bifurcation parameter, we prove the existence of local Hopf bifurcations about the zero equilibrium. For some critical parameter sets, the absolute local stability of the zero equilibrium is also guaranteed. The normal forms are calculated to determine the direction of the Hopf bifurcation and the stability of the bifurcating periodic solutions. By applying the theoretical results obtained, we prove the existence of bifurcating periodic solutions for a Nicholson's blowflies model and a Mackey-Glass model, both with a delayed harvesting term. We conclude that when the birth rate, death rate and capture rate are regulated, the maturation delay leads to Hopf-bifurcations as long as the difference between the maturation delay and the capture delay remains constant. The numerical simulations are presented illustrating the results. It is worth mentioning here that Theorems \ref{stability}, \ref{bifurcation} and \ref{estabilidaddireccion} can be applied to a wide class of biological models with two delays as long as the assumptions of these theorems are satisfied.

\section*{Acknowledgements}
J. Oyarce acknowledges support from Chilean National Agency for Research and Development (PhD. 2018--21180824) and 
A. G\'omez the support from the research  project 2120134 IF/R (University of B\'io-B\'io). 
The authors thank L. M. Villada for the advice in the numerical implementations of this work.

\section*{Declarations of Competing Interest} 
The authors declare they have no competing financial or personal interests that could  influence this paper. 

\section*{Funding}

This research did not receive any specific grant from funding agencies in the public, commercial or not-for-profits sectors. 

\section*{ORCID}

Adri\'an G\'omez \   https://orcid.org/0000-0002-2978-4465

Jos\'e Oyarce  \  https://orcid.org/0000-0002-0974-3463


\begin{thebibliography}{jose}

\bibitem{Amster} P. Amster, A. D\'eboli, Existence of $T$--periodic solutions of a generalized Nicholson's blowflies model with a nonlinear harvesting term,  \textit{App. Math. Lett.}  25(9) (2012) 1203-1207.

\bibitem{Qi} Q. An, E. Beretta, Y. Kuang, C. Wang, H. Wang,  Geometric stability switch criteria in delay differential equations with two delays and delay dependent parameters, \textit{J. Diff. Eqs.}266(11) (2019) 7073-7100.

\bibitem{Belair} J. B\'elair, S.A. Campbell, Stability and bifurcations of equilibria in a multiple-delayed differential equation,  \textit{SIAM J. Appl. Math.} 54(5) (1994) 1402-1424. 

\bibitem{berbrav2} L. Berezansky, E. Braverman, A note on stability of Mackey-Glass equations with two delays, \textit{J. Math. Anal. Appl.} 450(2) (2017) 1208-1228.

\bibitem{braverman} L. Berezansky, E. Braverman, L. Idels,  Nicholson's blowflies differential equations revisited: Main results and open problems, \textit{Appl. Math. Model.} 34 (2010) 1405-1417. 

\bibitem{bravd} E. Braverman, D. Kinzebulatov, Nicholson's blowflies equation with a distributed delay, \textit{Can. Appl. Math. Q.} 14(2) (2006) 107-128. 

\bibitem{Braddock} R.D. Braddock, P. van den Driessche, On a two lag differential delay equation, \textit{J. Austral. Math. Soc. Ser. B.}  24(3) (1983), 292-317.

\bibitem{ChowHale} S.N. Chow, J.K. Hale, Methods of bifurcation theory, Springer, New York, 1982.

\bibitem{Dieudonne} J. Dieudonn\'e, Foundations of Modern Analysis, Academic Press, New York, 1960.

\bibitem{Faria} T. Faria T, L.T. Magalh\~aes, Normal forms for retarded functional differential equations with parameters and applications to Hopf bifurcation, \textit{J. Diff. Eqs.} 122(2) (1995) 181-200.

\bibitem{Gopalsamy} K. Gopalsamy, Global stability in the delay-logistic equation with discrete delays, \textit{Houston J. Math.} 16 (1990) 347-356.

\bibitem{Gopalsamy3} K. Gopalsamy, Stability and oscillations in delay differential equations of population dynamics, Kluwer Academic Publishers,  Dordrecht, 1992.

\bibitem{Gopalsamy2} K. Gopalsamy, M.R.S Kulenovi\'c, G. Ladas,  Oscillations and global attractivity in models of hematopoiesis, \textit{J. Dyn. Diff. Eqs.} 2(2) (1990) 117-132.

\bibitem{Gu} K. Gu, S.L. Niculescu, J. Chen, On stability crossing curves for general systems with two delays, \textit{J. Math. Anal. Appl.} 311(11) (2005) 231-252.

\bibitem{Gurney} W.S.C Gurney, S.P. Blythe, R.M. Nisbet,  Nicholson's blowflies revisited, 
\textit{Nature.} 287 (1980) 17-21.

\bibitem{Gyori} I. Gy\"ori, G. Ladas, Oscillation theory of delay differential equations with applications, Oxford University Press, New York, 1991.

\bibitem{halegeometric} J.K. Hale, W. Huang, Global geometry of the stable regions for two delay differential equations,  \textit{J. Math. Anal. Appl.} 178 (1993) 344-362. 

\bibitem{Hale2} J.K. Hale, L.T. Magalh\~aes, W.M. Oliva, Dynamics in infinite dimensions, Second edition, Springer, New York, 2002.

\bibitem{Hale1} J.K. Hale, S.M. Verduyn Lunel, Introduction to functional differential equations, Springer, New York, 1993. 


\bibitem{Hassard} B.D. Hassard, N.D. Kazarinoff, Y.H. Wan, Theory and applications of Hopf bifurcation, Cambridge University Press, Cambridge, 1981.   


\bibitem{Huang2019} C. Huang, X. Yang, J. Cao, Stability analysis of Nicholson's blowflies equation with two different delays,  \textit{Math. Comput. Simul.} 171(9) (2020) 201-206. 

\bibitem{Lainscsek1} C. Lainscsek, L. Schettino, P. Rowat, E. van Erp, D. Song, H. Poizner, Nonlinear DDE analysis of repetitive hand movements in Parkinson's disease, pp. 421-427 in book Applications of Nonlinear Dynamics, Understanding Complex Systems (V. In, P. Longhini and A. Palacios eds.), Springer-Verlag,  Berlin Heidelberg, 2009.

\bibitem{Lainscsek2} C. Lainscsek, A.L. Sampson, R. Kim, M.L. Thomas, K. Man,  X. Lainscsek,  \emph{et. al}, Nonlinear dynamics underlying sensory processing dysfunction in schizophrenia,   \textit{Proc. Natl. Acad. Sci.}116 (2019) 3847-3852.

\bibitem{RuanWeiLi} X. Li, S. Ruan, J. Wei, Stability and bifurcation in delay-differential equations with two delays, \textit{J. Math. Anal. Appl.} 236(2) (1999) 254-280. 

\bibitem{Wei2021} Y. Liu, J. Wei, Bifurcation analysis in delayed Nicholson blowflies equation with delayed harvest, \textit{Nonlinear Dyn.} 105 (2021) 1805-1819.


\bibitem{Mackey1} M. Mackey, L. Glass, Oscillation and chaos in physiological control systems,  \textit{Science} 197 (1977) 287-289.


\bibitem{Nicholson} A.J. Nicholson, An outline of the dynamics of animal populations, \textit{Aus. J. Zool.} 2(1) (1954) 9-65. 


\bibitem{Piotrowska} M.J. Piotrowska, A remark on the ODE with two discrete delays,  \textit{J. Math. Anal. App.} 329(1) (2007) 664-676. 

\bibitem{ruan} S. Ruan, J. Wei, On the zeros of transcendental functions with applications to stability of delay differential equations with two delays, \textit{Dynam. Cont. Discr. Impul. Sys. Series A.} 10(6) (2003) 863-874.  


\bibitem{Hsmith} H. Smith, An introduction to delay differential equations with applications to the life sciences, Springer, New York, 2011. 

\bibitem{Wei2005} J. Wei, M.Y. Li, Hopf bifurcation analysis in a delayed Nicholson blowflies equation, \textit{Nonlinear Anal.} 60(7) (2005), 1351-1367.


\bibitem{Wei2007} J. Wei, D. Fan, Hopf bifurcation analysis in a Mackey-Glass system, \textit{Internat. J. Bifur. Chaos} 17 (2007) 2149-2157.



\end{thebibliography}
\end{document}